\RequirePackage{fix-cm}

\documentclass{svjour3}                    

\usepackage{psfrag,cite,amsmath,graphicx,epsfig,latexsym,subfigure,color}
\usepackage[scriptsize,bf]{caption}
\usepackage{color,xcolor}
\usepackage{amssymb}
\usepackage{setspace}
\usepackage{array}
\usepackage{booktabs}

\def\bn{\hfill \\ \smallskip\noindent}

\newcommand{\beq}{\begin{equation}}
\newcommand{\eeq}{\end{equation}}

\newcommand{\cI}{\mathcal{I}}

\newcommand{\trace}{{\mbox{\textrm{\rm Tr}}}}

\begin{document}
\def\pn {\par\smallskip\noindent}
\def \bn {\hfill \\ \smallskip\noindent}
\newcommand{\fs}{f_1,\ldots,f_s}
\newcommand{\f}{\vec{f}}
\newcommand{\hf}{\hat{f}}
\newcommand{\hx}{\hat{x}}
\newcommand{\hy}{\hat{y}}
\newcommand{\hz}{\hat{z}}
\newcommand{\hw}{\hat{w}}
\newcommand{\tw}{\tilde{w}}
\newcommand{\hlambda}{\hat{\lambda}}
\newcommand{\hbeta}{\hat{\beta}}
\newcommand{\tG}{\widetilde{G}}
\newcommand{\tg}{\widetilde{g}}
\newcommand{\barhx}{\bar{\hat{x}}}
\newcommand{\vecx}{x_1,\ldots,x_m}
\newcommand{\xoy}{x\rightarrow y}
\newcommand{\barx}{{\bar x}}
\newcommand{\bary}{{\bar y}}
\newcommand{\hrho}{\widehat{\rho}}
\newtheorem{thm}{Theorem}
\newtheorem{rmk}{Remark}[section]

\def\br{\break}
\def\smskip{\par\vskip 5 pt}
\def\proof{\bn {\bf Proof.} }
\def\QED{\hfill{\bf Q.E.D.}\smskip}
\def\qed{\quad{\bf q.e.d.}\smskip}
\newcommand{\cE}{\mathcal{E}}
\newcommand{\cM}{\mathcal{M}}
\newcommand{\cN}{\mathcal{N}}
\newcommand{\cJ}{\mathcal{J}}
\newcommand{\cT}{\mathcal{T}}
\newcommand{\bx}{\mathbf{x}}
\newcommand{\bp}{\mathbf{p}}
\newcommand{\bX}{\mathbf{X}}
\newcommand{\bY}{\mathbf{Y}}
\newcommand{\bP}{\mathbf{P}}
\newcommand{\bA}{\mathbf{A}}
\newcommand{\bB}{\mathbf{B}}
\newcommand{\bfM}{\mathbf{M}}
\newcommand{\bL}{\mathbf{L}}
\newcommand{\bz}{\mathbf{z}}
\newcommand{\cF}{\mathcal{F}}
\newcommand{\cR}{\mathcal{R}}
\newcommand{\bzero}{\mathbf{0}}

\newcommand{\blue}{\color{blue}}
\newcommand{\red}{\color{red}}

\title{A Linearly Convergent Doubly Stochastic Gauss-Seidel Algorithm for Solving Linear Equations and A Certain Class of  Over-Parameterized  Optimization Problems}
\author{Meisam Razaviyayn,$^*$ Mingyi Hong,$^*$ \\Navid Reyhanian, and Zhi-Quan Luo}
\date{Submitted in April 2018}
\institute{$^*$ equal contributions. M.\ Razaviyayn is with the Department of Industrial and Systems Engineering, the University of Southern California. \email{razaviya@usc.edu};
	M.\ Hong and N.\ Reyhanian are with the Department
	of Electrical and Computer Engineering, University of Minnesota, USA.
	\email{\{mhong, navid\}@umn.edu};
Z.-Q. Luo is with The Shenzhen Research Institute of Big Data, The Chinese University of Hong Kong, Shenzhen, China. \email{luozq@cuhk.edu.cn} This research is supported by the NSFC grants 61731018 and 61571384, the Peacock project of SRIBD.   
	}

 \newcommand{\bb}{\mathbf{b}}

	
	\maketitle
	
	\begin{abstract}
Consider the classical problem of solving a general linear system of equations $Ax=b$. It is well known that the (successively over relaxed) Gauss-Seidel scheme and many of its variants may not converge when $A$ is neither diagonally dominant nor symmetric positive definite. Can we have a linearly convergent G-S type algorithm that works for {\it any} $A$? In this paper we answer this question affirmatively by proposing a doubly stochastic G-S algorithm that is provably linearly convergent (in the mean square error sense) for any feasible linear system of equations. The key in the algorithm design is to introduce a {\it nonuniform  double stochastic} scheme for picking the equation and the variable in each update step as well as a stepsize rule. These techniques also generalize to certain iterative alternating projection algorithms for solving the linear feasibility problem $A x\le b$ with an arbitrary $A$, as well as  high-dimensional  minimization problems for training over-parameterized models in machine learning. Our results demonstrate that a carefully designed randomization scheme can make an otherwise divergent G-S algorithm converge.
		
	\end{abstract}
	
	\newpage
	\section{Introduction: Solving Linear System of Equations}
	Consider the generic problem of solving a linear system of equations
	\begin{align}\label{eq:linear:system}
	A x = b,	
	\end{align}
	where $A \in \mathbb{R}^{m\times n}$, $x, \; b \in \mathbb{R}^n$. We assume this system of equations has at least one solution. We use $[n]$ and $[m]$ to denote the set $\{1,\cdots, n\}$ and $\{1,\cdots, m\}$, respectively. A classical approach to solve \eqref{eq:linear:system} is by the Gauss-Seidel (G-S) algorithm, whereby at each iteration only {\it one} variable is updated by using the information from only {\it one} equation \cite{ortega1990numerical}. More precisely, let $A_{j:}$ and $b_j$ denote the $j$th row of $A$ and $j$th element of $b$, respectively. We define the G-S type algorithm as follows. 
	
	\begin{definition}\label{def:G-S}
	{\it An iterative algorithm for solving \eqref{eq:linear:system} is of {\it Gauss-Seidel type}, if at each iteration the algorithm updates one variable $x_i$, $i\in[n]$, by utilizing only $(A_{j:}, b_j)$ for some $j\in[m]$.}  
	\end{definition}

		A natural application of the G-S type algorithms is in the setting where $n$ players play a game in which $x_i$ is the variable of player~$i$. In this case we have $m=n$, and the objective of player~$i$ is to satisfy the $i$-th equation~$\sum_{j=1}^n a_{ij}x_j = b_i$. The best response strategy for player~$i$ is given by
		\begin{equation}
		\widehat{x}_i = \frac{b_i - \sum_{j\neq i}a_{ij}x_j}{a_{ii}}.
		\end{equation}
		Using a step-size $\alpha\ge 0$, this best response strategy leads to the following successive over-relaxation  (SOR) update rule:
		\begin{equation} \label{eq:BestResp}
		x_i^{r+1} = (1-\alpha) x_i^r + \alpha \frac{b_i - \sum_{j\neq i}a_{ij}x_j^r}{a_{ii}},
		\end{equation}
		where $a_{ij}$ is the $(i,j)$-th element of $A$; $b_i$ is the $i$th element of $b$.
		The central question is how to choose the step-size~$\alpha$ and determine the order in which the players update their variables so that the G-S type algorithm \eqref{eq:BestResp} will eventually lead to an equilibrium (or equivalently, a solution to \eqref{eq:linear:system}). 

	\subsection{Background on the Convergence of G-S Type Algorithm}

	To better understand the convergence behavior of G-S algorithm and its variants, let us consider the following example.\\
	
	\noindent\textbf{Example 1}: Consider the following $2\times 2$ special case of  \eqref{eq:linear:system}:
	\[
   A = \left[
	\begin{array}{cc}
	1 &-\tau\\
	-\tau &1\\
	\end{array}
	\right] \quad
	\textrm{and} \quad
	b = \left[
	\begin{array}{c}
	0  \\
	0\\
	\end{array}
	\right],
	\]
	where $\tau>1$ is some given constant. The best response strategy in \eqref{eq:BestResp} leads to the following update rule:
	\begin{align}
	x_1^{r+1} &= (1-\alpha) x_1^r + \alpha \tau x_2^r, \quad \mbox{when~$x_1$ is updated before $x^2$} \label{eq:x1:first}\\
	x_2^{r+1} &= (1-\alpha) x_2^r + \alpha \tau x_1^r, \quad \mbox{when~$x_2$ is updated before $x^1$.} \label{eq:x2:first}
	\end{align}
	Assume $x_1^0=x_2^0 > 0$.
	
	Let us consider the following five different update rules which are all of the G-S type. 
	\begin{enumerate}
	\item {\bf  Cyclic Successive Over Relaxation (SOR):} At iteration $r+1$, we perform
		\begin{align}
		x_1^{r+1} &= (1-\alpha) x_1^r + \alpha \tau x_2^r, \quad x_2^{r+1} = (1-\alpha) x_2^r + \alpha \tau x_1^{r+1}. 
		\end{align}
			\item {\bf  Symmetric SOR:} At iteration $r+1$, the variables are updated using a forward-sweep G-S step followed by a backward-sweep G-S step \cite{Saad03book}:
			\begin{align*}
			x_1^{r+1/2} &= (1-\alpha) x_1^r + \alpha \tau x_2^r, \quad x_2^{r+1/2} = (1-\alpha) x_2^r + \alpha \tau x_1^{r+1/2}. \\
			x_2^{r+1} &= (1-\alpha) x_2^{r+1/2} + \alpha \tau x_1^{r+1/2}, \quad x_1^{r+1} = (1-\alpha) x_1^{r+1/2} + \alpha \tau x_r^{r+1}. 
			\end{align*}
	\item {\bf Uniformly Randomized (UR) SOR:} At iteration $r+1$, randomly pick one variable from $\{x_1,x_2\}$ with equal probability. Update according to \eqref{eq:x1:first} or \eqref{eq:x2:first} based on which variables are selected, while fixing the other variable at its previous value.
	\item {\bf Non-Uniformly Randomized (NUR) SOR:} At iteration $r+1$, let $p^{r+1}_1>0$ and $p^{r+1}_2>0$ satisfy $p^{r+1}_1+p^{r+1}_2=1$;  randomly pick $x_i$ according to $p^{r+1}_i$. Update according to \eqref{eq:x1:first} or \eqref{eq:x2:first} based on which variables are selected, while fixing the remaining variable at its previous value.
	\item {\bf Random Permutation (RP) SOR:} At iteration $r+1$, randomly select a permutation $\pi$ of the index set $\{1,2\}$; The variables are updated according to
	\begin{align}\label{eq:radom:permute}
	x_{\pi(1)}^{r+1} &= (1-\alpha) x_{\pi(1)}^r + \alpha \tau x_{\pi(2)}^r, \quad x_{\pi(2)}^{r+1} = (1-\alpha) x_{\pi(2)}^r + \alpha \tau x_{\pi(1)}^{r+1}. 
	\end{align}
	Note that this method is referred to as the {\it shuffled} SOR in  \cite{Oswald15b}.
	\end{enumerate}
	It is easily seen that for any update order listed above, the resulting algorithm have the following property:
	\[
	\min\{x_1^r,x_2^r\} > \min\{x_1^0,x_2^0\}>0, \quad \forall r,\;\forall \alpha>0.
	\]
	On the other hand, the solution of the system of linear equation is $x_1^* = x_2^* = 0$;  hence none of these algorithms will find the solution. 
	\hfill $\blacksquare$
	
	\vspace{0.3cm}
	
	Next we give a brief literature review on the convergence analysis of the G-S type, as well as other related algorithms for solving a linear system of equations or  inequalities.
	
%

\noindent {\bf The convergence of G-S type algorithm.} It is well-known that when $A$ is either diagonally dominant, or symmetric positive definite (PD), then the classical SOR method with cyclic update rule converges to the solution of \eqref{eq:linear:system}; see \cite[Chapter 7]{ortega1990numerical}  and  \cite[Propositions 6.7, 6.8, 6.10]{bertsekas97}. More specifically, if $A$ is symmetric and PD, the convergence of SOR [for any $\alpha\in(0,2)$] can be established by showing that each iteration of the SOR algorithm is equivalent to a step of the coordinate descent (CD) algorithm for minimizing the strictly convex cost function $\frac{1}{2}x^TAx-b^Tx$ \cite[Section 2.6.3]{bertsekas97}. Note that the convergence rate in this case is linear, although the rate is not easily expressible in terms of the condition number of matrix $A$ if the classical cyclic rule is used \cite{Golub96}. Without the symmetry or the positive definiteness of $A$, the convergence of the SOR algorithm is only known when $A$ is diagonally dominant \cite[Section 2.6.2]{bertsekas97}. In particular,
using the matrix splitting $A=L+D+U$ where $L$, $U$, $D$ are the lower-triangular, upper-triangular and the diagonal part of $A$, respectively, we can write the SOR iteration as
\begin{align}
x^{r+1} = (1 + \alpha D^{-1} L )^{-1}\left[(1-\alpha)I -\alpha D^{-1} U\right] x^r + \alpha\left(I+\alpha D^{-1}L\right)^{-1}D^{-1}b.
\end{align}
Hence, the convergence of the SOR algorithm is guaranteed if the spectral norm of the iteration matrix is strictly less than one. Recently,  the work \cite{Oswald15b} shows that when $A$ is symmetric and positive semidefinite (PSD), then the G-S algorithm with RP rule can yield better convergence rate compared with the cyclic G-S (in the asymptotic region where $n$ is large).   From the above discussion it is clear that the classical G-S type algorithm does not work for any matrix $A$. A natural question is: Can a G-S type algorithm converge for any matrix $A$?

\vspace{0.2cm}

\noindent {\bf The convergence of Kaczmarz type algorithm.} Another popular method that bears similarity to the G-S type method for iteratively solving \eqref{eq:linear:system} is the Kaczmarz method \cite{KACZMARZ93}, whose iteration is expressed as
\begin{align}
x^{r+1} = x^r + \frac{ b_i - \langle A_{i:}, x^r\rangle }{\|A_{i:}\|^2} A^T_{i:}
\end{align}
where $A_{i:}$ is the $i$th row of $A$. This method has been used in many applications, but its rate of convergence was only analyzed in 2008 by Strohmer and Vershynin \cite{Strohmer2008} who proposed a randomized Kaczmarz (RK) method for over-determined linear systems. In the RK method, the $i$-th equation is selected for update randomly with probability proportional to $\|A_{i:}\|^2$. 
	This method can be seen as a particular case of stochastic gradient descent algorithm for minimizing the cost function
\[
\frac{1}{2}\frac{ (b_i - \langle A_{i:}, x\rangle)^2 }{\|A_{i:}\|^2},
\]
and the iterates converge linearly to a solution of \eqref{eq:linear:system}. The RK method has a convergence rate dependent only on a certain scaled condition number of matrix~$A$.
	
In a related work \cite{Leventhal10}, Leventhal and Lewis studied randomized variants of two classical algorithms, one is the CD for solving systems of linear equations (as  has been discussed above), the other is the iterated projection \cite{Deutsch01} for systems of linear inequalities (which contains the RK algorithm \cite{Strohmer2008} as a special case). The authors show that for the first algorithm when $A$ is symmetric (of size $n\times n$) and PSD, and for the second algorithm when the system has nonempty solution set, the global linear convergence can be established. Further the authors show that the linear rate can be  bounded in terms of natural linear-algebraic condition numbers of the problems. Note that for the iterated projection method  and the RK algorithm, the global linear convergence does not require $A$ to have full column rank.

\vspace{0.2cm}

Other recent works along this line include \cite{WHEATON2017, Gower15, DeLoera17,Needell10, Needell13, Oswald15}. In  \cite{Gower15},  a stochastic dual ascent (SDA) algorithm, which contains RK as a special case, was introduced for finding the projection of a given vector onto the solution space of a linear system. The method is dual in nature, with the dual being a unconstrained concave quadratic maximization problem. In each iteration of SDA, a dual variable is updated by choosing a point in a subspace spanned by the columns of a random matrix drawn independently from a fixed distribution. In \cite{DeLoera17}, the authors combined the relaxation
method of Motzkin \cite{Motzkin54} (also known as Kaczmarz method with the ``most violated constraint control") and the randomized Kaczmarz method \cite{Strohmer2008} to obtain a family of algorithms called {\it Sampling Kaczmarz-Motzkin} (SKM) for solving the linear systems $Ax\le b$. In SKM, at each time a subset of inequalities are picked, and the variables are updated based on the projection to the subspace corresponding to the most violated linear equality/inequality.
The reference \cite{WHEATON2017} proposed a new algorithm in which each iteration consists of obtaining $\ell_{\infty}$ norm projection of current approximate solution (onto hyperplanes defined by individual equations), followed by proper combination of the projections to all equations to yield the next iterate. Different from the Kaczmarz method \cite{KACZMARZ93}, this method requires information from all equations  to update one variable.  Needell \cite{Needell10} extended the RK method to the case of inconsistent equations, and showed that global linear convergence can be achieved until some fixed convergence horizon is reached. Needell and Tropp \cite{Needell13} analyzed a block version of the RK algorithm,  in which at each iteration the iterate is projected  onto the solution space of many equations simultaneously by selecting a block of rows rather than a single row. The convergence rate of the resulting algorithm is analyzed using the notion of {row paving} of a matrix.  Recently Liu and Wright proposed schemes to accelerate the RK method. The resulting scheme converges faster than the RK algorithm on ill conditioned problems \cite{Liu13ARK}.

Recently, there are a few works analyzing the randomized CD method  proposed in \cite{Leventhal10}  and the RK method \cite{Strohmer2008}, see \cite{Hefny17,MA15}. It is shown in \cite{MA15} that the randomized CD method can be extended to yield the minimum norm solution. In \cite{Hefny17}, variants of RK and randomized CD for solving Tikhonov regularized regression is proposed, and the corresponding convergence rates are derived. The rates derived indicate that RK based methods are preferable when $n > m$, while the randomized CD based methods are preferable when $m>n$.

It has been recognized that \emph{randomization} can be effective in simplifying the analysis of  RK method. 
In particular, Leventhal and Lewis \cite{Leventhal10} have used randomization in the RK method and strengthened the convergence analysis of the resulting randomized algorithm. They concluded that  ``randomization here provides a framework for simplifying the analysis of algorithms, allowing easy bounds on the rates of linear convergence in terms of natural linear-algebraic condition measures...". The present paper goes a step further in trying to understand the power of randomization. 
	\begin{center}
		\noindent\fcolorbox{black}[rgb]{0.9,0.9,0.9}{\begin{minipage}{1\textwidth}
				\begin{center}
					{\bf (Q1)} ~~Can randomization make the (otherwise divergent) G-S algorithm convergent for a general linear system \eqref{eq:linear:system}?
				\end{center}
			\end{minipage}}
		\end{center}

\subsection{Contribution of This Work}
	In this paper, we answer the above question affirmatively. In particular, we propose a \emph{doubly stochastic G-S algorithm} that is provably linearly convergent (in expectation) for any feasible linear system of equations. The key in the algorithm design is to introduce a {\it nonuniform  double randomization} scheme for picking the equation and the variable in each update step of the G-S algorithm, along with an appropriate stepsize rule. Interestingly, these randomization techniques also generalize to certain iterative alternating projection algorithms for solving the linear feasibility problem $A x\le b$ with an arbitrary $A$, as well as to certain high-dimensional over-parameterized minimization problems. Our results demonstrate that a carefully designed randomization scheme can make an otherwise divergent G-S algorithm converge linearly.

	\subsection{Notation}
	For any matrix $A$, let $\|A^{\dagger}\|_2$ denote the smallest constant $M$ such that
	$\|A x\|_2\ge \frac{1}{M}\|x\|_2$ for all $x$. Let us define the {\it relative condition number} of $A$ as $k(A):=\|A\|_2 \|A^{\dagger}\|_2$; the {\it scaled condition number} is defined as $\kappa(A):=\|A\|_F\|A^{\dagger}\|_2$. It is easy to verify that
	\begin{align}
	1\le \frac{\kappa (A)}{\sqrt{n}}\le k(A). \label{eq:kappakrelation}
	\end{align}
	We use $A_{:i}$ and $A_{j:}$ to denote the $i$th column and $j$th row of $A$, respectively.  For a symmetric matrix $B$, we use $\lambda_{\max}(B)$, $\lambda_{\min}(B)$ and $\underline{\lambda}_{\min}(B)$ to denote its maximum, the minimum and the minimum nonzero eigenvalues.

%
	\section{A Doubly Stochastic G-S Algorithm}\label{sec:introduction}

	\subsection{Simultaneously selecting equations and variables?}
\begin{figure}[htb]
			\centering
			\includegraphics[width=9cm]{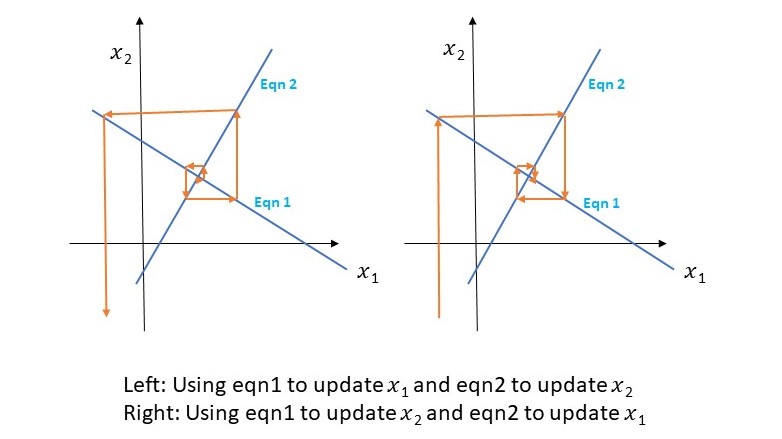}
			\caption{Convergence behavior of the G-S algorithm under two different variable to equation associations.} \label{fig:G-S}
		\end{figure}

Recall that in the G-S algorithm \eqref{eq:BestResp}, we always use equation $i$ to update variable $x_i$. In other words, variable $x_i$ is \emph{locked} to equation $i$. This locking is arbitrary since one can simply re-order the equations (or re-indexing the variables) without affecting the solution of the linear system, and yet different variable to equation coupling will give rise to a different G-S update scheme, some of which may be divergent while others may be convergent. Figure \ref{fig:G-S} gives an illustrative example in $\mathbb{R}^2$,  whereby if we use equation 1 to update variable $x_1$, and equation 2 to update variable $x_2$, then the G-S algorithm diverges (left subfigure), but if we use  equation 1 to update variable $x_2$, and equation 2 to update variable $x_1$, then the G-S algorithm converges linearly (right subfigure). This example suggests variable to equation association can greatly affect the convergence of the G-S algorithm.

Thus, a natural way to design a convergent G-S algorithm for a $n\times n$ general linear system \eqref{eq:linear:system} is to carefully select a fixed matching that determines which variable is to be updated by which equation. Clearly, the choice of a good matching (one that can lead to a convergent G-S algorithm) will be dependent on the coefficient matrix $A$. Unfortunately, this is a challenging task since the number of possible matchings is $n!$, which grows exponentially in $n$. Moreover, for a non-square linear system ($m\times n$), it is not clear how to define such a matching between variables and equations.

In this paper, we propose to \emph{unlock} the fixed pairing of each variable to a unique equation in the G-S algorithm. Moreover, since determining which variable is to be updated by which equation is hard,  we propose to simply do so \emph{randomly}!
More specifically, at each G-S iteration, we can {\it randomly} select a pair $(i,j)$ where $i$ is an index for an equation while $j$ is an index for a variable. Using this strategy, after picking the pair $(i,j)$, one can update variable~$x_j$ using equation $i$ as follows (according to the G-S update rule) 
	\begin{equation} \label{eq:BestRespij}
	x_j^{r+1} = (1-\alpha) x_j^r + \alpha \frac{b_i - \sum_{k\neq j}a_{ik}x_k^r}{a_{ij}}, \quad x^{r+1}_\ell = x^{r}_{\ell}, \; \forall~\ell\ne j.
	\end{equation}
 To answer $\rm (Q1)$, it is then natural to ask the following question:
		\begin{center}
			\noindent\fcolorbox{black}[rgb]{0.9,0.9,0.9}{\begin{minipage}{1\textwidth}
					\begin{center}
						{\bf (Q2)} ~~Can {\it unlocking} the fixed variable-equation pairing and \emph{randomization} ensure the convergence of a G-S type algorithm for an arbitrary linear system~\eqref{eq:linear:system}?
					\end{center}
				\end{minipage}}
				\vspace{0.1cm}
			\end{center}

	To understand the impact of unlocking and randomization, let us consider the following example.
	
		\noindent\textbf{Example 2}: Consider the same $(A,b)$ as given in Example 1, and let us consider the {\it unlocked} version of the G-S outlined above.
		
		Specifically, after selecting the pair $(i,j)$, we can use one of the following four update rules to update the variables:
	\begin{enumerate}
		\item {\bf Case 1.} $i=1,j=1\rightarrow$ $x_1^{r+1} = (1-\alpha) x_1^r + \alpha \tau x_2^r$;
		\item {\bf Case 2.} $i=1,j=2\rightarrow$ $x_2^{r+1} = (1-\alpha) x_2^r + \frac{\alpha}{\tau} x_1^r$;
		\item {\bf Case 3.} $i=2,j=1\rightarrow$ $x_1^{r+1} = (1-\alpha) x_1^r + \frac{\alpha}{\tau} x_2^r$;
		\item {\bf Case 4.} $i=1,j=2\rightarrow$ $x_2^{r+1} = (1-\alpha) x_2^r +\alpha\tau  x_1^r$.
	\end{enumerate}
	Consider a \emph{uniform randomized update rule} where at each iteration,  one of above update rules is selected (each with probability 1/4) and
used to update the variable $x$. Consider an initialization that $x_1$ and $x_2$ have the same sign. Without loss of generality
let us assume $x_1,x_2>0$.  Define the random process $z^r \triangleq \min\{x_1^r,x_2^r\}$. Using the uniform randomized update
rule, one can show that
\begin{equation}
z^{r+1} \geq \left\{
\begin{array}{lll}
(1-\alpha) z^r + \alpha \tau z^r &\textrm{with probability 1/4} & \textrm{scenario 1}\\
(1-\alpha) z^r + \frac{\alpha}{ \tau} z^r &\textrm{with probability 1/4} & \textrm{scenario 2}\\
z^r  &\textrm{with probability 1/2} & \textrm{scenario 3}\\
\end{array}
\right..
\end{equation}
In order to show the divergence of the algorithm, it suffices to show that $z^r \rightarrow \infty$ with probability one.
To show this, we define the random process $\{w^r\}_{r=0}^{\infty}$ with $w^0 = z^0$ and
\begin{equation}
w^{r+1} = \left\{
\begin{array}{ll}
(1-\alpha + \alpha \tau) w^r & \textrm{if scenario 1 happens in process }z^r\\
(1-\alpha + \frac{\alpha}{ \tau}) w^r  &  \textrm{if scenario 2 happens in process }z^r\\
w^r   &  \textrm{if scenario 3 happens in process }z^r\\
\end{array}
\right.
\end{equation}
Clearly, $z^r \geq w^r, \;\forall r$. Hence we only need to show $w^r \rightarrow \infty$ with probability one. Notice that $\log(w^r) = \sum_{i=1}^r \beta^i+ \log(w^0)$ where $\beta^i$ is an i.i.d process with
\begin{equation}
\beta^{i} \geq \left\{
\begin{array}{ll}
\log(1-\alpha + \alpha \tau)  & \textrm{with probability 1/4}\\
\log(1-\alpha + \frac{\alpha}{ \tau})  &  \textrm{with probability 1/4}\\
0   &  \textrm{with probability 1/2}\\
\end{array}
\right.
\end{equation}
It is not hard to see that $\mathbb{E}[\beta^i]>0, \;\forall i$. Therefore, $\lim_{r\rightarrow \infty} \sum_{i=1}^r \beta^i = \infty$  due to the law of large numbers. Consequently, $\lim_{r\rightarrow} \log(w^r) = \infty$ which implies  $\lim_{r\rightarrow \infty} \|x^r\| = \infty$, regardless of stepsize $0<\alpha<1$.

Furthermore, if one uses the cyclic update rule, then at each iteration of the G-S algorithm, each one of the above cases will be selected once according to some deterministic rules. Therefore, in order to study the convergence of the resulting algorithm, one need to look at the spectral radius of the resulting mapping. For example, if we select the update rules in the order of 1), 2), 3), and 4), one needs to study the spectral radius  $\rho(B_4B_3B_2B_1)$ where
\[
B_1 = \left[
\begin{array}{cc}
1-\alpha \tau & \alpha \tau \\
0 & 1\\
\end{array}
\right],  \;\;
B_2 = \left[
\begin{array}{cc}
1 & 0 \\
\alpha/\tau & 1-\alpha\\
\end{array}
\right],\;\;
B_3 = \left[
\begin{array}{cc}
1-\alpha  & \alpha/\tau \\
0 & 1\\
\end{array}
\right],\;\;
B_4 = \left[
\begin{array}{cc}
1& 0 \\
\alpha \tau & 1-\alpha\\
\end{array}
\right].
\]
One can check that any permutation of the above matrices will result in a product matrix whose spectral radius is larger than 1. Consequently, the resulting algorithm will diverge for almost all initialization, and for any $0<\alpha<1$. Furthermore, one can numerically check that the randomly permuted rule will also diverge for almost all initialization and for any $0<\alpha<1$.
%
	\hfill $\blacksquare$\\

\begin{figure}[htb]
			\centering
			\includegraphics[width=9cm]{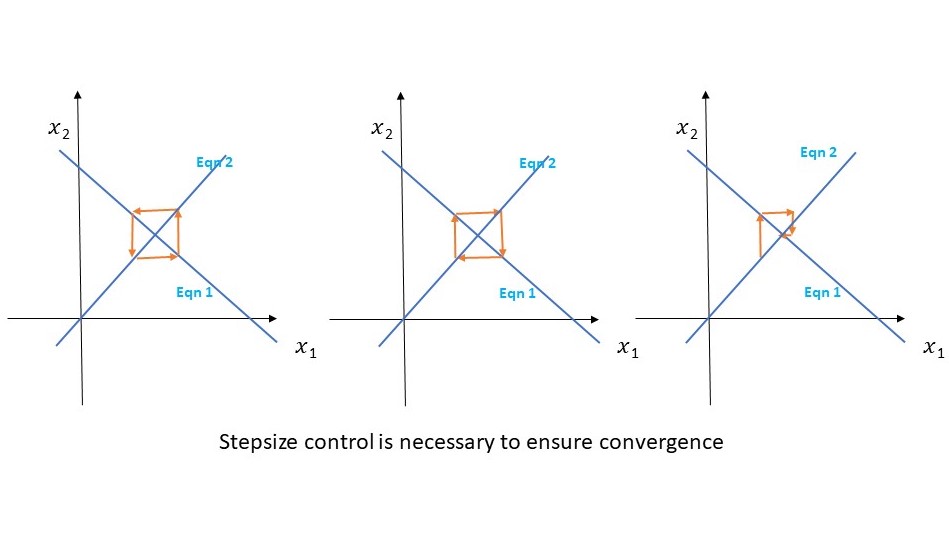}
			\caption{Convergence behavior of the G-S algorithm under different stepsizes} \label{fig:G-S_stepsize}
		\end{figure}

	The above example suggests that by simply unlocking the variable-equation association, the G-S type methods still may not converge. Moreover, Figure~\ref{fig:G-S_stepsize} shows an example in $\mathbb{R}^2$ whereby the G-S type algorithm will not converge if $\alpha=1$ regardless of variable-equation association, but will converge if $0<\alpha<1$ for any variable-equation association. Thus, stepsize control is also necessary (in addition to randomization) for the convergence of the G-S type algorithm.

Surprisingly, we will show in the subsequent sections that, by properly selecting a data-dependent updating probability as well as the stepsize $\alpha$, the randomized unlocked G-S algorithm (Algorithm 1 below) is globally linearly convergent in the mean squared error sense for any feasible linear system \eqref{eq:linear:system}. In this algorithm, at each iteration only one index pair $(i,j)$ is selected. However, the selection is  based on a non-uniform distribution defined over  index pairs.
	
\subsection{A Doubly Stochastic G-S Type Algorithm}
We propose a doubly stochastic G-S type algorithm below.  This randomized algorithm combines the {\it unlocking} idea with certain {\it non-uniform} random selection of both equations and variables. The key that ensures the convergence of the resulting algorithm, compared with the divergent cases in Example 2, is a judicious selection of the probability $p_{ij}$ that governs how the equations and the variables are picked, as well as a stepsize rule.  Note that if a particular entry of the matrix $a_{ij}=0$, then its corresponding $p_{ij}=0$.  Hence the corresponding $(i,j)$-th index will never be picked with probability one. Therefore the update \eqref{eq:x:update} is well defined.

	\begin{center}
		\fbox{
			\begin{minipage}{\linewidth}
				\smallskip
				\centerline{\bf {Algorithm 1. A doubly stochastic G-S (DSGS) algorithm.}}
				\smallskip Let $\alpha>0$. At iteration $0$, randomly generate $x^0$.
				
				At iteration $r+1$,  randomly select one index pair $(i,j)$ with the probability
				\begin{align}\label{eq:p}
				p_{ij} = \frac{a^2_{ij}}{\sum_{i',j'}a^2_{i'j'}}.
				\end{align}
				Update $x_{j}$ by the following:
				\begin{align}\label{eq:x:update}
				x^{r+1}_j = (1-\alpha) x_j^r + \alpha \left( \frac{b_i-\sum_{k\ne j} a_{ik} x_k^r}{a_{ij}}\right) = x_j^r+\alpha\left(\frac{b_i - \sum_{k=1}^{n} a_{ik}x_k^r}{a_{ij}} \right),
				\end{align}
				and for all $k$, $k\neq j$, set 
				\begin{align}
				x^{r+1}_k = x^r_k.
				\end{align}

			\end{minipage}
		}
	\end{center}
	
	
	\subsection{Case 1: $A$  has full column  rank}
	We first consider the case where $A$ has full column rank. Let $x^*$  be the feasible solution for \eqref{eq:linear:system}.
	For simplicity of notation we use the superscript $``+"$ to denote the new iteration. Let us define
	\begin{align}
	\Delta^+_j = x^{+}_j- x^*_j, \quad \Delta_j = x_j - x^*_j.
	\end{align}
	We have the following result.
	\begin{thm}\label{claim:linear:full:rank}
		Consider a consistent system $Ax=b$ with $A$ being full column rank. Then there holds
\begin{equation*}
  \mathbb{E}\left[\|\Delta^{+}\|^2\mid x\right] \le\Delta^T \left( I + \frac{n\alpha^2-2\alpha}{\sum_{i,j}a^2_{ij}} A^T A\right)\Delta.
\end{equation*}
Thus the DSGS algorithm converges globally linearly in the mean squared error sense for $0<\alpha<2/n$. Moreover, if we choose $\alpha=1/n$, then the DSGS algorithm achieves the following convergence rate:
		\begin{align}\label{eq:linear:fullrank}
		\mathbb{E}\left[\|\Delta^{+}\|^2\mid x\right] \le \left(1-\frac{1}{n \kappa^2(A)}\right)\|\Delta\|^2 \le \left(1-\frac{1}{n^2  k^2(A)}\right)\|\Delta\|^2.
		\end{align}
	\end{thm}
	\noindent{\bf Proof.} Clearly we  have the following relation
	\begin{align}\label{eq:opt}
	a_{ij} x^*_j = b_i -\sum_{k\ne j}a_{ik}x^*_k, \quad \forall~i.
	\end{align}
	Also let us choose $p_{ij} = \frac{a^2_{ij}}{\sum_{i,j}a^2_{ij}}, \; \forall\; i,j. $	
	We have the following series of equalities
	\begin{align*}
	&\mathbb{E}[(\Delta^{+}_j)^2\mid x] \nonumber\\
	&= \left(1-\sum_{i}p_{ij}\right)(x_j-x^*_j)^2 +\sum_{i}p_{ij}\left((1-\alpha) x_j +\alpha\frac{b_i-\sum_{k\ne j}a_{ik}x_k}{a_{ij}} -x^*_j\right)^2\nonumber\\
	&= \left(1- \frac{\sum_{i}a^2_{ij}}{\sum_{i,j}a^2_{ij}}\right)\Delta^2_j +\frac{1}{\sum_{i,j}a^2_{ij}}\sum_{i}\left(a_{ij}(1-\alpha) x_j + \alpha(b_i -\sum_{k\ne j}a_{ik}x_k) -
	(1-\alpha)x^*_j a_{ij} -\alpha x^*_j a_{ij}\right)^2\nonumber\\
	&\stackrel{\eqref{eq:opt}}=  \left(1- \frac{\sum_{i}a^2_{ij}}{\sum_{i,j}a^2_{ij}}\right)\Delta^2_j +\frac{1}{\sum_{i,j}a^2_{ij}}\sum_{i}\left(a_{ij}(1-\alpha) \Delta_j + \alpha(b_i -\sum_{k\ne j}a_{ik}x_k) -
	\alpha(b_i - \sum_{k\ne j} a_{ik} x^*_k)\right)^2\nonumber\\
	&=\left(1- \frac{\sum_{i}a^2_{ij}}{\sum_{i,j}a^2_{ij}}\right)\Delta^2_j +\frac{1}{\sum_{i,j}a^2_{ij}}\sum_{i}\left(a_{ij} \Delta_j-\alpha \sum_{k}a_{ik}\Delta_k\right)^2\nonumber\\
	&=\left(1- \frac{\sum_{i}a^2_{ij}}{\sum_{i,j}a^2_{ij}}\right)\Delta^2_j +\frac{1}{\sum_{i,j}a^2_{ij}}\sum_{i}\left(a^2_{ij} \Delta^2_j + \alpha^2 (\sum_{k}a_{ik}\Delta_k)^2 - 2\alpha\sum_{k}a_{ik}\Delta_k a_{ij}\Delta_j\right)\nonumber\\
	& = \Delta^2_j + \frac{1}{\sum_{i,j}a^2_{ij}}\left(\sum_{i}\alpha^2 (\sum_{k}a_{ik}\Delta_k)^2 - 2 \alpha \Delta_j\sum_{i}\sum_{k}a_{ik}\Delta_k a_{ij}\right)\nonumber\\
	& = \Delta^2_j +  \frac{1}{\sum_{i,j}a^2_{ij}}\left(\alpha^2\sum_{i}(A_{i:} \Delta)^2 - 2 \alpha \Delta_j A^T_{:j} A \Delta\right).
	\end{align*}
	In the last equation, the notation $A_{i:}$ denotes the $i$th row of the matrix $A$.
	
	Summing the above equation over $j$, we have
	\begin{align}
	\mathbb{E}\left[\|\Delta^+\|^2\mid x\right] &= \|\Delta\|^2 + \alpha^2 \frac{n}{\sum_{i,j}a^2_{ij}} \Delta^T A^T A \Delta - \frac{2\alpha}{\sum_{i,j}a^2_{ij}} \sum_{j} \Delta_j A^T_{:,j}A \Delta\nonumber\\
	& = \|\Delta\|^2+ \alpha^2 \frac{n}{\sum_{i,j}a^2_{ij}} \|\Delta\|^2_{A^T A} - \frac{2\alpha}{\sum_{i,j}a^2_{ij}}  \|\Delta\|^2_{A^T A}\nonumber\\
	&= \Delta^T \left( I + \alpha^2 \frac{n}{\sum_{i,j}a^2_{ij}}  A^T A - \frac{2\alpha}{\sum_{i,j}a^2_{ij}}  A^T A\right)\Delta\nonumber\\
	&= \Delta^T \left( I + \left(\alpha^2 \frac{n}{\sum_{i,j}a^2_{ij}}   - \frac{2\alpha}{\sum_{i,j}a^2_{ij}} \right) A^T A\right)\Delta\nonumber.
	\end{align}
	Clearly, to make the error converge geometrically, it suffices to have
	\begin{align*}
	\alpha^2 n  - 2\alpha<0, \quad \mbox{or equivalently}\quad \alpha <\frac{2}{n}.
	\end{align*}
	Let us pick $\alpha = 1/n$, then we have
	\begin{align*}
	\mathbb{E}[\|\Delta^+\|^2\mid x] &\le  \|\Delta\|^2 \left( 1  - \frac{1}{n\sum_{i,j}a^2_{ij}} \lambda_{\min}(A^T A)\right) \nonumber\\
	&=  \|\Delta\|^2 \left( 1 - \frac{1}{n \|A\|^2_F \|A^{\dagger}\|_2^2}\right) = \|\Delta\|^2 \left( 1  - \frac{1}{n \kappa^2(A)} \right) .
	\end{align*}
	This shows that the expected value of the optimality gap shrinks globally geometrically.  \QED
	
	\begin{remark}
		Let us compare the rate obtained above with existing results in the literature. First, the rate of the randomized CD method obtained in \cite[Theorem 3.4]{Leventhal10}  (for solving \eqref{eq:linear:system} with $A$ being symmetric and PD) is given by

		\begin{align}\label{eq:pd:rate}
		 \left(1-\frac{1}{\|A^{\dagger}\|_2\trace[A]}\right) 
		&\le \left(1-\frac{1}{\sqrt{n}\kappa (A)}\right)\le \left(1-\frac{1}{n k  (A)}\right),
		\end{align}
		where the last inequality is due to~\eqref{eq:kappakrelation}.
		We can see that our rate obtained in \eqref{eq:linear:fullrank} takes a similar form, except that our rate is proportional to $\left(1-\left(\frac{1}{\sqrt{n}\kappa (A)}\right)^2\right)$. This is  due to the lack of symmetry and positive definiteness of $A$.\\
		
		Alternatively, the rate obtained in \cite{Strohmer2008} when using RK method for solving \eqref{eq:linear:system} with $A$ being full column rank is given by (see, e.g., \cite[Theorem 4.2]{Leventhal10}, \cite[Proposition 2]{DeLoera17})
		\begin{align}\label{eq:pd:rate:rk}
		\mathbb{E}[\|\Delta^{+}\|^2\mid x]\le \left(1-\frac{2\alpha-\alpha^2}{\kappa^2(A)}\right)\|\Delta\|^2 \stackrel{\alpha = 1}= \left(1-\frac{1}{\kappa^2(A)}\right)\|\Delta\|^2.
		\end{align}
		Note that at each iteration of $RK$, $n$ variables are updated. In contrast, our rate in \eqref{eq:linear:fullrank} is proportional to $\left(1-\frac{1}{ n \kappa^2(A)}\right)$, but at each iteration only one variable is updated. When $n$ is large and $\kappa(A)$ is large, we have
		\begin{align*}
		\left(1-\frac{1}{ n \kappa^2(A)}\right)^n \approx \exp(-1/\kappa^2(A)) \approx  1- \frac{1}{\kappa^2(A)}
		\end{align*}
		which indicates that the two rates are asymptotically comparable.  \hfill $\blacksquare$
		\end{remark}
	
	\subsection{Case II: $A$ has no full column rank}
	In this subsection, we assume that $A$ has no full column rank, and system \eqref{eq:linear:system} has a least one solution $x^*$. In this case, the previous analysis does not work because \eqref{eq:linear:fullrank} does not imply linear convergence as $\kappa(A) = \infty$. In this section, we use a different analysis technique.
	
	Let us define
	\begin{align}\label{eq:beta}
	\beta: = A x -b \in\mathbb{R}^{m}, \quad \mbox{with}\quad \beta_k = A_{k:} x -b_k.
	\end{align}
	Below we will show that the quantity  $\|A (x-x^*)\|^2 = \|\beta\|^2$ converges linearly to zero in expectation.
		\begin{thm}\label{claim:linear:not:full:rank}
			Consider a consistent system $Ax=b$ with arbitrary $A$. Let us pick
			$$\alpha=  \frac{1}{\|A\|^2_F}\underline{\lambda}_{\min}(A A^T),$$
			and define $\beta^+ := Ax^+ - b$  to be the residual after one update of the algorithm.
			Then the double stochastic G-S algorithm achieves the following convergence rate
	\begin{align}\label{eq:linear:rank:deficit}
	\mathbb{E}[\|\beta^{+}\|^2 \mid x] & \le \|\beta\|^2 \left(1-  \left(\frac{1}{\|A\|^2_F}\underline{\lambda}_{\min}(A A^T)\right)^2\right).
	\end{align}
		\end{thm}
		
	\noindent{\bf Proof.} First from the update rule \eqref{eq:x:update}, it is clear that when the tuple $(i,j)$ is picked, we have
	\begin{align}
	x^{+} = x -\alpha \left( \frac{A_{i:} x - b_i}{a_{ij}} \right) e_j
	\end{align}
	where $e_j$ is the $j$th elementary vector. Left multiplying both sides by $A_{k:}$, we have
	\begin{align}
	A_{k:}x^{+} = A_{k:} x - \alpha \left( \frac{A_{i:} x -b_i}{a_{ij}}\right) a_{kj}.
	\end{align}
	According to the definition \eqref{eq:beta}, we further have
	\begin{align}
	\beta^{+}_k &= A_{k:}(x^{+}-x^*) = (A_{k:}x^{+}-b_k)\nonumber\\
	& = \beta_k - \alpha\frac{A_{i:}x-b_i}{a_{ij}}a_{kj}= \beta_k - \alpha\frac{\beta_i}{a_{ij}}a_{kj}.
	\end{align}
Since each tuple $(i,j)$ is picked using probability $p_{ij}$ in \eqref{eq:p}, we have the following estimate
	\begin{align}
	\mathbb{E}[\|\beta^{+}_k\|^2\, | \, x]& = \sum_{i,j}p_{ij} \left(\beta_k-\alpha\frac{\beta_i}{a_{ij}}a_{kj}\right)^2= \frac{1}{\sum_{i,j}a^2_{ij}}\sum_{i,j}\left(\beta_k a_{ij}-\alpha{\beta_i}{a_{kj}}\right)^2.
	\end{align}
	Summing over $k$, we have
	\begin{align}
	\mathbb{E}[\|\beta^{+}\|^2\, | \, x] &= \frac{1}{\sum_{i,j}a^2_{ij}}\sum_{k,i,j}\left(\beta_k a_{ij}-\alpha{\beta_i}{a_{kj}}\right)^2\nonumber\\
	&= \frac{1}{\sum_{i,j}a^2_{ij}}\sum_{k,i,j}\left(\beta^2_k a^2_{ij}-2 \beta_k a_{ij}\alpha{\beta_i}{a_{kj}} + \alpha^2 a^2_{kj} \beta^2_i\right)\nonumber\\
	& = \frac{1}{\sum_{i,j}a^2_{ij}}\left((1+\alpha^2) \|A\|^2_F\|\beta\|^2 -2 \alpha  \sum_{k,i,j} (\beta_i a_{ij}) (\beta_k {a_{kj}})\right)\nonumber\\
	& = \frac{1}{\sum_{i,j}a^2_{ij}}\left((1+\alpha^2) \|A\|^2_F\|\beta\|^2 -2 \alpha  \beta^T A A^T \beta \right)\nonumber\\
	& = \beta^T \left( (1+\alpha^2)I - \frac{2\alpha }{\sum_{i,j}a^2_{ij}} A A^T \right)\beta \nonumber.
	\end{align}
	{Note that by definition $\beta:= A x -b$, and the system is consistent. Therefore $\beta$ belongs to the column space of $A A^T$}. It follows that we have
	$$\|A^T \beta\|^2\ge \|\beta\|^2\underline{\lambda}_{\min}(AA^T).$$
	Therefore, the following sufficient conditions are needed
	\begin{align}
	(1+\alpha^2) - 2\frac{\alpha}{\sum_{i,j}a^2_{ij}}{\lambda}_{\max}(A A^T)> 0\nonumber\\
	(1+\alpha^2) - 2\frac{\alpha}{\sum_{i,j}a^2_{ij}}\underline{\lambda}_{\min}(A A^T)<1.
	\end{align}
By picking the value
	\begin{align}\label{eq:alpha}
	\alpha=  \frac{1}{\sum_{i,j}a^2_{ij}}\underline{\lambda}_{\min}(A A^T),
	\end{align}
	 we have
	\begin{align*}
	\mathbb{E}[\|\beta^{+}\|^2 \, | \, x] & \le \|\beta\|^2 \left(1-  \left(\frac{1}{\sum_{i,j}a^2_{ij}}\underline{\lambda}_{\min}(A A^T)\right)^2\right). 
	\end{align*}
	The claim is proved. \hfill \QED

	\begin{remark}\label{rmk:non:fullrank}
		Let us compare the rate obtained above with the rate of a randomized CD method (Algorithm 3.5 in \cite{Leventhal10}), a method which updates only one variable at each iteration, while utilizing one column of matrix $A$.  It is shown that for a consistent linear systems of equations \eqref{eq:linear:system} with arbitrary non-zero matrix $A$, the randomized CD method achieves the  rate 
%
	
	\begin{align}\label{eq:rate:psd}
	\mathbb{E}[\|\beta^+\|_{A^TA}^2\mid x]\le \left(1-\frac{\underline{\lambda}_{\min}(A^T A)}{\|A\|^2_F}\right)\|\beta\|_{A^TA}^2.
	\end{align}
	Clearly the above rate is closely related to the one given in \eqref{eq:linear:rank:deficit}. \hfill $\blacksquare$
\end{remark}

{ \begin{remark}
	Although we are mainly interested in addressing Question {\bf Q1}, namely how to design a G-S scheme that converges for any matrix A, the proposed double stochastic algorithm does have important practical value. Consider the distributed computation setting where there is a central controller  node $0$ connected to a number of distributed computing nodes, each having a subset of rows of data matrix $A$ and $b$, respectively. The central controller stores the variable $x$, the error term $Ax -b$ and is capable of broadcasting to every distributed nodes. At each iteration, the central node randomly pick a pair $(i,j)$, and send the distributed node that has $A_{i:}$ the scalar $A_{i:} x -b_i$; the corresponding node will update $x_j$ according to \eqref{eq:x:update} using its {\it local} information. Then the new $x_i$ will be transmitted back to node $0$, and node $0$ will recompute $A x-b$ and continue the previous process. Therefore after $n$ iterations of the algorithm, the total number of messages transmitted between the local nodes to node $0$ is $2\times n$. In comparison, if one implements the RK method in the same distributed network, then each iteration $n$ messages have to be communicated from the local nodes to node $0$.
\end{remark}}

	\section{A Doubly Stochastic Alternating Projection Algorithm}
	In this section, we extend our previous analysis to the problem of finding a point in the intersection of multiple polyhedral sets. The algorithm to be developed  has the flavor of the classical alternating projection algorithm, except that we perform the alternating projection coordinate-wise.
	Specifically, we consider the following problem:
	\begin{align}\label{eq:projection}
	\mbox{Find}~x\quad \quad \mbox{s.t.}\quad A_{i:} x\le b_i, \quad i=1,\cdots, m,
	\end{align}
where we write the system of inequalities $Ax \leq b$ in the above form to emphasize the existence of $m$ separate scalar inequalities. We also assume that  this system of linear inequalities is feasible.

\vspace{0.2cm}

Our proposed algorithm is closely related to Algorithm 1, except that we only update those inequalities that are violated.
	\begin{center}
		\fbox{
			\begin{minipage}{\linewidth}
				\smallskip
				\centerline{\bf {Algorithm 2. The double stochastic alternating projection algorithm.}}
				\smallskip At iteration $0$, randomly generate $x^0$.
				
				At iteration $r+1$, randomly pick an index pair $(i,j)$ with probability
				$$p_{ij} = \frac{a^2_{ij}}{\sum_{i,j}a^2_{ij}}.$$
				
				Update $x_{j}$ by the following
				\begin{align}\label{eq:x:update:2}
				x^{r+1}_j &= x_j^r, \quad \mbox{if}~A_{i:} x^r\le b_i\\
				x^{r+1}_j &= (1-\alpha) x_j ^r+ \alpha \left( \frac{b_i-\sum_{k\ne j} a_{ik} x_k^r}{a_{ij}}\right) \nonumber\\
				&= x_j^r+\alpha\left(\frac{b_i - \sum_{k=1} a_{ik}x_k^r}{a_{ij}} \right), \quad \mbox{otherwise}.
				\end{align}
				For all $k, k \neq j$, set $x_k^{r+1} = x_k^r$.
				
			\end{minipage}
		}
	\end{center}
	
	To facilitate our analysis, let us define the following function
	\begin{align}\label{eq:obj}
	f(x): =\sum_{i=1}^{m}f_i(x),
	\end{align}
	where
	\[
	f_i(x) :=  \frac{1}{2}(A_{i:} x - b_i)^2_{+}.
	\]
	We note that any feasible solution of \eqref{eq:projection} will imply $f(x)=0$. Further, each function $f_i$ is  differentiable, and its gradient is $A^T_{i:}(A_{i:} x-b_i)$ if $A_{i:} x - b_i\ge 0$, and it is $0$ otherwise. 
	For a given iteration $r$, define the index set
	\begin{align}
	\mathcal{I}^r:=\{i\mid A_{i:} x^r - b_i> 0\},
	\end{align}
	and define $\Omega^r_{\cI}\in\mathbb{R}^{m\times m}$ as the diagonal matrix with $\Omega^r_{\cI}[i,i]=1$ if $i\in\mathcal{I}^r$ and $\Omega^r_{\cI}[i,i]=0$ otherwise.
	Then we have
	\begin{align}
	\sum_{i=1}^{m}f_i(x^r) = \sum_{i\in{\cI}^r}f_i(x^r) = \frac{1}{2}\|A^r x^r -b^r\|^2
	\end{align}
	where $A^r= \Omega_{\cI}^r A \in\mathbb{R}^{m\times n}$, and $b^r$ is defined similarly. 
	
	Note that by using the above definition, we have
	\begin{align}\label{eq:key:1}
	\sum_{i=1}^{m}f_i(x^r)=  \frac{1}{2}\|A^r x^r -b^r\|^2 \le  \frac{1}{2}\|A x^r - b\|^2.
	\end{align}

	 \subsection{Case 1: $A$ has full row rank}
	In this subsection we  make the following assumption
	\begin{align}
	\lambda_{\min}(A A^T)>0.
	\end{align}
	

	Similarly as before, we will use $x^{+}$ (resp. $x$) to denote the new (resp. previous) iteration; 
	we will use $A_\cI$, $b_{\cI}$ and $\Omega_{\cI}$ to denote $A^r, b^r, \Omega_{\cI}^r$ at iteration $r$, respectively.
	
	\begin{thm}\label{claim:full:rank:projection}
		Suppose $A$ has full row rank, and $\alpha$ is chosen as
		\begin{align}
		\alpha <\frac{\lambda_{\min}(A A^T)}{\|A\|_F^2}.
		\end{align}
		Then we have
		\begin{align}
		\mathbb{E}[f(x^{+})\mid x] \le \left(1- \left(\frac{\lambda_{\min}(A A^T)}{2 \|A\|^2_F}\right)^2\right) f(x)
		\end{align}
	\end{thm}
	
	\noindent{\bf Proof.} Suppose that the $(i,j)$th pair gets selected, and that $j$th coordinate gets updated (this means that $i\in\cI$), then we can rewrite $x^+$ as following
	\begin{align}
	x^+ = x+\alpha\left(\frac{b_i - \sum_{k} a_{ik}x_k}{a_{ij}} \right) e_j.
	\end{align}
	In this case, we can estimate the component function $f_{\ell}(x^{+})$ based on whether the $\ell$th inequality is satisfied for $x$. Suppose that $\ell\in \cI$ (i.e. the $\ell$-th inequality is not satisfied), then we have
	\begin{align}
	f_\ell(x^+) &= f_\ell\left(x+\alpha\left(\frac{b_i - \sum_{k} a_{ik}x_k}{a_{ij}} \right) e_j\right)\nonumber\\
	&\le f_\ell(x) + \left\langle \nabla f_\ell(x), \alpha\left(\frac{b_i - \sum_{k} a_{ik}x_k}{a_{ij}} \right)e_j\right\rangle + \frac{1}{2}\left(\alpha a_{\ell j}\frac{b_i - \sum_{k} a_{ik}x_k}{a_{ij}} \right)^2\nonumber\\
	&= f_\ell(x) + \alpha \left\langle A^T_{\ell:}(A_{\ell:} x -b_\ell), \left(\frac{b_i - \sum_{k} a_{ik}x_k}{a_{ij}} \right)e_j\right\rangle + \frac{1}{2}\left(\alpha a_{\ell j}\frac{b_i - \sum_{k} a_{ik}x_k}{a_{ij}} \right)^2\nonumber\\
	&= f_\ell(x) + \alpha \left\langle a_{\ell j}(\sum_{k}a_{\ell k} x_k -b_\ell) e_j, \left(\frac{b_i - \sum_{k} a_{ik}x_k}{a_{ij}} \right)e_j\right\rangle \nonumber\\
	&\quad + \frac{\alpha^2}{2}\left(a_{\ell j}\frac{b_i - \sum_{k} a_{ik}x_k}{a_{ij}} \right)^2\label{eq:mean:1st}.
	\end{align}
	Otherwise, if $\ell\notin \cI$ (i.e. the $\ell$-th inequality is satisfied), we have
	\begin{align}
	f_\ell(x^+) &= f_\ell\left(x+\alpha\left(\frac{b_i - \sum_{k} a_{ik}x_k}{a_{ij}} \right) e_j\right)\nonumber\\
	& = \frac{1}{2}\left(A_{\ell:} x -b_\ell +\alpha A_{\ell:}\left(\frac{b_i - \sum_{k} a_{ik}x_k}{a_{ij}} \right) e_j\right)_+^2\nonumber\\
	&\le f_\ell(x)  + \frac{1}{2}\left(\alpha a_{\ell j}\frac{b_i - \sum_{k} a_{ik}x_k}{a_{ij}} \right)^2
	\label{eq:mean:2nd},
	\end{align}
	where the last inequality is due to the fact that for all $\ell\notin\cI$, $A_{\ell :}x -b_\ell\le 0$. Further, if $(i,j)$-th pair gets selected, but that $j$-th coordinate is {\it not} updated, then we have $x^+ = x$.

	Using the above inequalities, we have  the following 
	\begin{align*}
	&\mathbb{E}\left[\sum_{\ell}f_\ell(x^{+}) \mid x\right] =  \sum_{\ell\in \cI}\mathbb{E}\left[ f_\ell(x^{+})\mid x \right] +  \sum_{\ell\notin \cI}\mathbb{E}\left[ f_\ell(x^{+})\mid x \right]  \nonumber\\
	&\stackrel{\rm{(i)}}\le \sum_{\ell\in\cI} \sum_{ij: i\in\cI}p_{ij} \left( \alpha \left\langle a_{\ell j}\left(\sum_{k}a_{\ell k} x_k -b_\ell\right) e_j, \left(\frac{b_i - \sum_{k} a_{ik}x_k}{a_{ij}} \right)e_j\right\rangle\right) \nonumber\\
	&\quad + \sum_{\ell}\sum_{ij: i\in\cI}p_{ij}\frac{\alpha^2}{2}\left(a_{\ell j}\frac{b_i - \sum_{k} a_{ik}x_k}{a_{ij}} \right)^2+ \sum_{\ell}\sum_{i,j} p_{ij}f_\ell(x)\\
	&= \frac{\alpha}{\sum_{i,j}a^2_{ij}}\sum_{\ell\in\cI} \sum_{ij: i\in\cI} \left\langle a_{\ell j}\left(\sum_{k}a_{\ell k} x_k -b_\ell\right) e_j, a_{ij}\left({b_i - \sum_{k} a_{ik}x_k} \right)e_j\right\rangle \nonumber\\
	&\quad + \sum_{\ell}\sum_{ij:i\in\cI}\frac{\alpha^2}{2\sum_{i,j}a^2_{ij}}\left(a_{\ell j} (b_i - \sum_{k} a_{ik}x_k) \right)^2 + \sum_{\ell}f_\ell(x)\\
%
& \stackrel{\rm{(ii)}}\le \sum_{\ell\in\cI}f_\ell(x) -\frac{\alpha}{\sum_{i,j}a_{ij}^2} \|A^T_{\cI}\left(A_{\cI} x -b_{\cI}\right)\|^2 + \frac{\alpha^2}{2}\|A_{\cI}x -b_{\cI}\|^2\\
	& \le \left(1- \frac{\alpha\lambda_{\min}(A_{\cI} A^T_{\cI})}{\|A\|_F^2} + \frac{\alpha^2}{2}\right)\sum_{\ell}f_\ell(x) \le \left(1- \frac{\alpha\lambda_{\min}(A A^T)}{\|A\|_F^2} + \frac{\alpha^2}{2}\right)\sum_{\ell}f_\ell(x)
	\end{align*}
where in $\rm (i)$ we have used the two cases \eqref{eq:mean:1st} and \eqref{eq:mean:2nd}; in $\rm(ii)$ we have used the fact that $f_{\ell}(x)= 0$ for $\ell\not\in {\cal I}$; in the last inequality we have used the  fact that $A_{\cI} A^T_{\cI}$ is a principal submatrix of $A A^T$, the fact that $f_\ell(x)\ge 0$, and the fact that $\alpha$ is chosen small enough such that
	$$-\frac{\alpha\lambda_{\min}(A A^T )}{\|A\|_F^2} + \frac{\alpha^2}{2}<0.$$
	This concludes the proof. \QED
	
	\subsection{Case 2: $A$ does not have full row rank}
	In this subsection, we present an analysis of Algorithm 2 without the full row rankness assumption. To this end, we need to use the well-known Hoffman's error bound.
	\begin{lemma}
		Let $\mathcal{S}$ denote the solution set for the linear system in the constraint \eqref{eq:projection}. Then there exists a constant $\tau>0$ independent of $b$, with the following property
		\begin{align}\label{eq:hoffman}
		x\in\mathbb{R}^n, \; S\ne \emptyset \rightarrow \; \mbox{dist}(x, S) \le \tau \| (Ax-b)^{+}\|,
 		\end{align}
 		where we have defined
 		\begin{align}
 		(Ax-b)_i^{+} =\max\{0, A_{i:} x - b\}, \quad \mbox{dist}(x,S):=\inf_{y\in S}\|x-y\|.
 		\end{align}

	\end{lemma}

	Assume that the system \eqref{eq:projection} is feasible, and let $S$ denote its solution set and let $x^*\in S$. Clearly we have $f(x^*):=\sum_{i=1}^{m}f_i(x^*)=0$.  We have the following claim.
	\begin{thm}
		Consider a feasible system $Ax\le b$ with arbitrary $A$. Let us pick
		$\alpha=  \frac{1}{n}.$
		Then Algorithm 2 achieves the following convergence rate
		\begin{align}\label{eq:linear:rank:deficit:inequality}
		\mathbb{E}[\mbox{dist}^2(x^{+}, S)\mid x] \le \left(1- \frac{1}{ n\tau^2\sum_{i,j}a^2_{ij}} \right) \mbox{dist}^2(x, S).
		\end{align}
	\end{thm}
	
	\noindent{\bf Proof.} Recall that from the update rule we have
	\begin{align}
	x^+ = x+\alpha\left(\frac{b_i - \sum_{k} a_{ik}x_k}{a_{ij}} \right) e_j.
	\end{align}
	Let us define the projection of $x$ to the feasible set as
	\begin{align}\label{eq:projection:x}
	P(x):= \arg\min_{y\in S}\|x-y\|.
	\end{align}
	We have the following relationship for $\mathbb{E}\left[\mbox{dist}(x^{+}, S)^2\mid x\right]$
	\begin{align*}
	&\mathbb{E}[\mbox{dist}^2(x^{+}, S)\mid x] \nonumber\\
	& = \mathbb{E}\left[\|x^{+}-P(x^+)\|^2\mid x\right] \le \mathbb{E}\left[\|x^{+}-P(x)\|^2\mid x\right]\\
	&= \mathbb{E}\left[\|x^{+}-x\|^2\mid x\right] + 2 \mathbb{E}\left[\langle x^{+}-x, x-P(x)\rangle\mid x\right] +\|x-P(x)\|^2 	\label{eq:exp:error}
	\end{align*}
	where the inequality is due to the definition of the projection \eqref{eq:projection:x}.
	Let us bound the above equality term by term. First we have
	\begin{align}
	\mathbb{E}\left[\|x^{+}-x\|^2\mid x\right]& = \sum_{(i,j):i\in\mathcal{I}} p_{ij} \alpha^2 \left\|\left(\frac{b_i - \sum_{k} a_{ik}x_k}{a_{ij}} \right) e_j\right\|^2\nonumber\\
	&= \sum_{(i,j):i\in\mathcal{I}} \frac{1}{\sum_{i,j}a^2_{ij}} \alpha^2 \left\|\left({b_i - \sum_{k} a_{ik}x_k} \right) e_j\right\|^2\nonumber\\
	&= \sum_{i\in\mathcal{I}} \frac{1}{\sum_{i,j}a^2_{ij}} \alpha^2 n \left\|{b_i - \sum_{k} a_{ik}x_k} \right\|^2\nonumber\\
	&= \frac{\alpha^2 n}{\sum_{i,j}a^2_{ij}} \sum_{i\in\mathcal{I}}  \left\|A_{\mathcal{I}}x - b_{\mathcal{I}}\right\|^2\nonumber.
	\end{align}
	The second term in \eqref{eq:exp:error} is given by
	\begin{align}
	\mathbb{E}\left[\langle x^{+}-x, x-P(x)\rangle \mid x \right] &= 2 \alpha \sum_{(i,j): i\in\mathcal{I}}\left\langle p_{ij} \left(\frac{b_i - \sum_{k} a_{ik}x_k}{a_{ij}} \right) e_j, x-P(x)\right\rangle \nonumber\\
	&= 2 \alpha \frac{1}{\sum_{i,j}a^2_{ij}}\sum_{(i,j): i\in\mathcal{I}}\left\langle \left(a_{ij}(b_i - \sum_{k} a_{ik}x_k)\right) e_j, x-P(x)\right\rangle \nonumber\\
	&= 2 \alpha \frac{1}{\sum_{i,j}a^2_{ij}}\left\langle A^T_{\mathcal{I}} (b_{\mathcal{I}}-A_{\mathcal{I}}x), x-P(x)\right\rangle \nonumber\\
	&= - 2 \alpha \frac{1}{\sum_{i,j}a^2_{ij}}\left\langle \nabla f(x), x-P(x)\right\rangle \nonumber\\
	&\le - 2 \alpha \frac{1}{\sum_{i,j}a^2_{ij}}(f(x)-f(P(x))) \nonumber\\
	&= - 2 \alpha \frac{1}{\sum_{i,j}a^2_{ij}}\|A_{\mathcal{I}}x -b_{\mathcal{I}}\|^2, \nonumber
	\end{align}
	where the first equality is due to the fact that $x-P(x)$ is constant when conditioned on $x$; the first inequality is due to the convexity of $f$; and the last equality is because the system is feasible so $f(P(x))=0$.
	
	Therefore, overall we have
	\begin{align}
	\mathbb{E}[\mbox{dist}^2(x^{+}, S)\mid x] \le \frac{n \alpha^2 - 2 \alpha}{\sum_{i,j}a^2_{ij}}\|A_{\mathcal{I}}x -b_{\mathcal{I}}\|^2 + \|x-P(x)\|^2\nonumber.
	\end{align}
	Therefore if $0<\alpha<  2 /n$, we can apply the Hoffman condition \eqref{eq:hoffman}
		\begin{align}
		\mathbb{E}[\mbox{dist}^2(x^{+}, S) \mid x] &\le \frac{n \alpha^2 - 2 \alpha}{\tau^2\sum_{i,j}a^2_{ij}} \mbox{dist}^2(x, S) + \|x-P(x)\|^2\nonumber\\
		& = \left(1 + \frac{n \alpha^2 - 2 \alpha}{\tau^2\sum_{i,j}a^2_{ij}} \right) \mbox{dist}^2(x, S) \nonumber\\
		& \stackrel{\alpha=\frac{1}{n}} \le \left(1- \frac{1}{ n\tau^2\sum_{i,j}a^2_{ij}} \right) \mbox{dist}^2(x, S).
		\end{align}
Therefore we conclude that the algorithm converges linearly in expectation. 	\QED
\begin{remark}
	Let us compare the above result with the rate for a randomized iterative projection method (Algorithm 4.6) developed in \cite{Leventhal10} (cf.  \cite[Theorem 4.7]{Leventhal10}), which is given below
			\begin{align}
			\mathbb{E}[\mbox{dist}^2(x^{+}, S)\mid s] \le \left(1- \frac{1}{ \tau^2\sum_{i,j}a^2_{ij}} \right) \mbox{dist}^2(x, S).
			\end{align}
	Note that the randomized iterative projection method updates $n$ variables at each iteration, therefore we can use an argument similar to  Remark \ref{rmk:non:fullrank} to see the two methods have comparable rates. \hfill $\blacksquare$
\end{remark}

	\section{Generalization to Double Stochastic Gradient Descent in Over-Parameterized Machine Learning }
	In this section, we extend the previous analysis and algorithm design to solve problems beyond linear system of equalities and inequalities.
	
	\vspace{0.4cm }
	
	Consider the following problem
	\begin{align}\label{eq:general:formulation}
	\min_{x}\; f(x):= \sum_{i=1}^{m}f_i(x_1,\cdots, x_n)
	\end{align}
	where each $f_i: \mathbb{R}^{n}\to \mathbb{R}$. We assume the following throughout this section.
	
	\noindent{\bf Assumption 1.} 
	\begin{enumerate}
		\item[(1)] Each function $f_i$ is non-negative, i.e.,
		$$f_i(x)\ge 0, \quad \forall~x\in\mathbb{R}^n.$$
		\item[(2)] Assume that each $f_i$ has Lipschitz continuous gradient with constant $L_i$, i.e.,
		\begin{align}\label{eq:lip:property}
		\|\nabla f_i(x) - \nabla f_i(z)\|\le L_i \|x-z\|, \quad \forall~x,z,\in\mathbb{R}^{n},\quad i = 1, \cdots, n.
		\end{align}
		Further assume that, for all $x\in \mathbb{R}^n$, $f$ satisfies
\begin{align}\label{eq:strong:convexity}
f(x)-f(P(x))\le \langle \nabla f(x), x-P(x)\rangle -\frac{\gamma}{2}\|P(x)-x\|^2, 
\end{align}
		where $P(x)$ is the projection of $x$ to the set of global minimizers for problem \eqref{eq:general:formulation} similarly as defined in \eqref{eq:projection:x}; $\gamma>0$ is some constant.
		\item[(3)] Assume that the global optimal objective value of problem \eqref{eq:general:formulation} is zero, i.e.,
		\begin{equation}
		\label{eq:MinZero}
		\min_{x} f(x)=0.
		\end{equation}
	\end{enumerate}
	\vspace{0.3cm}
	
	The assumption in \eqref{eq:strong:convexity} is slightly weaker than the usual strong convexity assumption.  For example, the remark below shows that the function defined in \eqref{eq:obj} satisfies all the conditions in Assumption~1. Additionally, for any optimal solution~$x^*$  for problem \eqref{eq:general:formulation}, Assumption~1 implies that
	\begin{align}\label{eq:gradient:zero}
	\nabla f_i(x^*) =0, \quad\forall~i.  
	\end{align}

\begin{remark}
{The function defined in \eqref{eq:obj} satisfies Assumption~1. To see this, one can first easily  verify that the function defined in \eqref{eq:obj}	satisfies  condition \eqref{eq:lip:property} with $L_i = \|a_i\|^2$. Moreover, for a given point $x$, define $\mathcal{I} = \{i\, | \, a_i^T x -b_i \geq 0\}$.  Noticing that $A_{\mathcal{I}}x - b_{\mathcal{I}} \geq 0$ and $A_{\mathcal{I}}P(x) - b_{\mathcal{I}} \leq 0$, we obtain
	\[
	\langle A_{\mathcal{I}}x - b_{\mathcal{I}} ,  A_{\mathcal{I}} P(x) - b_{\mathcal{I}} \rangle \leq 0.
	\]
	Adding the term $\|A_{\mathcal{I}}x - b_{\mathcal{I}}\|^2$ to both sides and rearranging the terms, we obtain
	\begin{equation}\label{eq:xxstarIneq}
	\|A_{\mathcal{I}}x - b_{\mathcal{I}}\|^2 \leq \langle   A_{\mathcal{I}}^T\left(A_{\mathcal{I}} x- b_{\mathcal{I}}\right) , x- P(x) \rangle.
	\end{equation}
	On the other hand, Hoffman error bound implies that there exists a constant $\gamma$ (independent of $x$) such that
	\begin{equation} \label{eq:HoffmanErrbd}
	\frac{\gamma}{2} \|x - P(x)\|^2 \leq \|A_{\mathcal{I}} x - b_{\mathcal{I}}\|^2.
	\end{equation}
	Multiplying \eqref{eq:xxstarIneq} by two and adding to \eqref{eq:HoffmanErrbd}, we get
	
	\[
	\|A_{\mathcal{I}} x - b_{\mathcal{I}}\|^2  +   \frac{\gamma}{2} \|x-P(x)\|^2\leq 2 \langle A_{\mathcal{I}}^T(A_{\mathcal{I}}  x - b_{\mathcal{I}})  , x -  P(x)\rangle,
	\]
	i.e., \eqref{eq:strong:convexity} holds. } \hfill $\blacksquare$
\end{remark}

\begin{remark}
We also note that the condition \eqref{eq:strong:convexity} is weaker than the traditional strong convexity, and it is also weaker than the essentially strong convexity condition defined in \cite{liu14asynchronous,karimi2016linear}.  In particular, the essentially strong convexity requires that 
	\begin{align}\label{eq:strong:convexity:essential}
	f(x)-f(P(x))\le \langle \nabla f(x), x-y\rangle -\frac{\gamma}{2}\|x-y\|^2, \; \forall~x, y,\;\mbox{s.t.}\quad P(x)= P(y).
	\end{align}
	The above condition clearly implies \eqref{eq:strong:convexity}. Moreover, it is not hard to see that for any strongly convex function $h(\cdot)$, and any linear mapping $A$, the function $f(x) = h(Ax)$ satisfies  \eqref{eq:strong:convexity}.   Notice that the condition~\eqref{eq:strong:convexity} does not even imply the convexity of the function~\cite{necoara2018linear}. This condition is referred to as \textit{quasi strong convexity} in \cite{necoara2018linear} and as \textit{weak strong convexity} in \cite{karimi2016linear}.\hfill $\blacksquare$	
\end{remark}

\begin{remark}
[\textbf{Connection to  over-parameterized  learning}] In training machine learning models via empirical risk minimization framework, the goal is to minimize the empirical risk loss
\begin{equation}
\label{eq:Overparameterized}
\min_{x} \;\; \sum_{i=1}^m\; \ell (a_i^T x ,y_i),
\end{equation}
where $\{(a_i,y_i)\}_{i=1}^m$ is the data and $x$ is the parameter of the model. The loss function $\ell(\cdot,\cdot)$ is non-negative and $\ell(y',y) = 0 $ if $y=y'$. In the over-parameterized regime where the number of parameters~$n$ is lager than the number of data points~$m$, it is not hard to see that the minimum value of \eqref{eq:Overparameterized} is zero and hence \eqref{eq:MinZero} is satisfied. \hfill $\blacksquare$
\end{remark}

%
%

	
Under the conditions of Assumption~1, let us consider the following gradient based algorithm algorithm:
	\begin{center}
		\fbox{
			\begin{minipage}{\linewidth}
				\smallskip
				\centerline{\bf {Algorithm 3. A doubly stochastic gradient algorithm.}}
				\smallskip At iteration $0$, randomly generate $x^0$.
				
				At iteration $r+1$, pick the index pair $(i,j)$ with probability  $p_{ij}=p=\frac{1}{mn}$.
				
				Update $x_{j}$ by the following
				\begin{align}\label{eq:x:update:3}
				x^{r+1}_j &= x_j^r - \alpha \nabla_j f_i(x^r),
				\end{align}
				and for all $k, k\neq j$, set
				\begin{align}
				x_k^{r+1} = x_k^r
				\end{align}
			\end{minipage}
		}
	\end{center}
In contrast to the previous two algorithms, now a uniform sampling probability is used in Algorithm 3. It turns out that for stochastic gradient based algorithm such a choice is sufficient to guarantee convergence.
	
		\begin{thm}
			Suppose Assumption 1 is satisfied. Then applying Algorithm 3 to problem \eqref{eq:general:formulation} with a stepsize 
\[
0< \alpha \le \frac{\gamma^2}{2 \sum_{\ell}L_{\ell} \sum_{i=1}^{n}L^2_i}
\]
achieves the following convergence rate
			\begin{align}\label{eq:general:rate}
			\mathbb{E}[f(x^{+})\mid x ] \le \left(1-\frac{\alpha \gamma}{mn}\right)  f(x).
			\end{align}
		\end{thm}
		\noindent{\bf Proof.} Suppose that $(i,j)$th pair gets selected, then we have
		\begin{align}
		x^{+} = x - \alpha \nabla f_i(x) e_j.
		\end{align}
		For a fixed $\ell\in[m]$, and by using \eqref{eq:lip:property} we have the following
		\begin{align}
		f_{\ell}(x^+) &= f_{\ell}(x-\alpha \nabla f_i(x)e_j)\nonumber\\
		&\le f_{\ell}(x) -\langle \nabla f_{\ell}(x), \alpha \nabla f_i(x)e_j\rangle +\frac{L_{\ell}\alpha^2}{2}\| \nabla f_i(x)e_j\|^2\nonumber.
		\end{align}
		Therefore we obtain the following relations:
		\begin{align}
		\mathbb{E}[f(x^{+})\mid x] &\le \frac{1}{mn}\sum_{\ell\in[m]}\sum_{i,j} \left( -\alpha \langle \nabla f_{\ell}(x) e_j, \nabla f_i(x) e_j\rangle + \frac{L_{\ell}\alpha^2}{2}\| \nabla f_i(x)e_j\|^2\right) + f(x)\nonumber\\
		& = \frac{1}{mn}\sum_{\ell\in[m]}\sum_{i,j} \left( -\alpha \langle \nabla f_{\ell}(x) e_j, \nabla f_i(x) e_j\rangle\right)  + \frac{ \sum_{\ell} L_{\ell}\alpha^2}{2 mn} \sum_{i=1}^{n}\| \nabla f_i(x)\|^2+ f(x)\nonumber\\
		& = f(x) - \frac{\alpha}{mn}\|\sum_{i=1}^{m}\nabla f_i(x)\|^2  + \frac{ \sum_{\ell}L_{\ell}\alpha^2}{2 mn} \sum_{i=1}^{n}\| \nabla f_i(x)\|^2\nonumber\\
		& \stackrel{\eqref{eq:gradient:zero}}= f(x) - \frac{\alpha}{mn}\|\nabla f(x)\|^2  + \frac{ \sum_{\ell}L_{\ell}\alpha^2}{2 mn} \sum_{i=1}^{n}\| \nabla f_i(x)-\nabla f_i(x^*)\|^2\nonumber\\
		& \stackrel{\eqref{eq:lip:property}}\le  f(x) - \frac{\alpha}{mn}\|\nabla f(x)\|^2  + \frac{ \sum_{\ell}L_{\ell}\alpha^2}{2 mn} \sum_{i=1}^{n} L^2_i \|x-x^*\|^2.\label{eq:intermmediate}
		\end{align}
		Using the property \eqref{eq:strong:convexity}, we have
		\begin{align}
		f(x)-f(x^*)& \le \langle \nabla f(x), x-x^*\rangle  -\frac{\gamma}{2}\|x-x^*\|^2 \nonumber\\
		&\le \frac{1}{\gamma} \|\nabla f(x)\|^2+ \frac{\gamma}{4}\|x-x^*\|^2-\frac{\gamma}{2}\|x-x^*\|^2 \nonumber\\
		&= \frac{1}{\gamma} \|\nabla f(x)\|^2- \frac{\gamma}{4}\|x-x^*\|^2.
		\end{align}
		Combine the above with the assumption that $f(x^*)=0$, we obtain
		\begin{align}
		-\frac{\alpha}{mn}\|\nabla f(x)\|^2 \le -  \frac{\alpha \gamma}{mn} f(x) -\frac{\gamma^2}{4}\frac{\alpha}{mn} \|x-x^*\|^2.
		\end{align}
		Plugging this inequality into \eqref{eq:intermmediate}, we obtain
		\begin{align}
		\mathbb{E}[f(x^{+})\mid x] \le \left(1-\frac{\alpha \gamma}{mn}\right) f(x) -\left(\frac{\gamma^2}{4}\frac{\alpha}{mn} - \frac{ \sum_{\ell}L_{\ell}\alpha^2}{2 mn} \sum_{i=1}^{N} L^2_i\right) \|x-x^*\|^2.
		\end{align}
		Therefore by choosing the following constant
		\begin{align}
		\alpha = \frac{\gamma^2}{2 \sum_{\ell}L_{\ell} \sum_{i=1}^{N}L^2_i}
		\end{align}
		we obtain the desired result. 		\QED

It should be noted that  stochastic block coordinate descent methods or primal-dual methods have been proposed in the literature for optimization problems where only one function and one variable is sampled in each iteration. However, in order to obtain a linear rate of convergence, they either need to explicitly involve the dual variables \cite{yu2015doubly}  and/or apply variance reduction techniques by maintaining a running average of the gradients \cite{zhang2016accelerated}. In contrast, Algorithm~3  works with only primal variables, does not maintain a running average of gradients, and does not require strong convexity assumption. In addition, in the case of system of linear equations, Algorithm~1 allows exact minimization per-iteration by non-uniform sampling of the equations and variables, while the works mentioned above does not necessarily converge with the exact optimization per iteration.
		

		\section{Numerical Results}
		In this section, we evaluate the performance of the doubly stochastic G-S algorithm proposed in this paper against some of the existing algorithms, for solving linear systems  of equations~\eqref{eq:linear:system}.  In particular, the performance of DSGS is compared with  S2CD  \cite{Konecny14},  RK \cite{Strohmer2008},  RGS   \cite{MA15}, and the REK method  \cite{zouzias2013randomized}.  Note that except for the S2CD algorithm, the rest of the algorithms are designed specifically for solving linear system of equations. The S2CD algorithm is  an extension of the SVRG algorithm \cite{Johnson13} to include coordinate descent update, and it solves a strongly convex finite sum problem using a doubly stochastic update. For fair comparison, in all algorithms, we increase the iteration counters when $n$ updates of the variable coordinates are completed. We divide our experiments into two broad categories: when $A^T\:A$ is PD and when $A^T\:A$ is PSD. 
			\begin{table}
			\centering
			\caption{\textsc{The Average Number of Required Iterations  to Reach the Error Level of $10^{-10}$ For An Overdetermined System with a Random Coefficient Matrix $A$.} }
			\begin{tabular}{llllll}
				\hline \hline
				Size ($m\times n$)& DSGS       & RK & RGS & REK     &  S2CD    \\
				\hline
				\hline
				$40\times 20$ & {\bf 3,612} & 4,334 & 4,342 & 4,949 & 8,671 \\
				$40\times 30$   & {\bf 17,211} & 29,653 & 30,174 & 34,067  & 59,892  \\
				$60\times 30$   & 7,715  & {\bf 7,205} & 7,230 & 8,348  & 14,392  \\
				$60\times 50$  & {\bf 36,503} & 41,638 & 41,829 & 47,545 & 81,542\\
				$80\times 40$  &  13,262 &  {\bf 9,026} & 9,058  &  10,591 & 18,195\\
				$80\times 60$  & {\bf 83,857} & 84,974 & 84,871 & 95,807 & 166,224\\
				\hline
			\end{tabular}
			\label{tab:res1}
			\vspace{1ex}
		\end{table}
		
		In the first set of experiments, we consider over-determined systems of equations with $m\ge n$. To implement S2CD, we consider an objective function $f(x)=\frac{1}{m}\sum_{i=1}^{m}\|A_{i:}x-b_i\|^2$ with $f_i(x)=\|A_{i:}x-b_i\|^2$. Similar to \cite{Strohmer2008,zouzias2013randomized}, the elements of matrix $A$ are independently drawn from standard Gaussian distribution~$N(0,1)$. Clearly, the condition number of such a matrix changes with the size. A randomly generated vector whose elements are drawn from standard Gaussian $N(0,1)$ is considered as $x^*$.  Accordingly, the vector $b$ is found such that $b=Ax^*$. We terminate an algorithm when it reaches an error value of $\|x-x^*\|\le 10^{-10}$. Each entry in the table is computed by averaging the results for running a given algorithm over ten randomly generated problems.  The result of this experiment is summarized in Table~\ref{tab:res1}.
As can be seen from this table, when matrix $A$ gets close to a square matrix, the number of required iterations increases. This is reasonable since the condition numbers of our random matrix increases.  Moreover, for these harder cases, $40\times 30$, $60\times 50$ and $80\times 60$, the proposed approach in this paper slightly outperforms other algorithms.

	We also test our algorithm on problems instances generated with different given condition numbers. To do so, we first generate a zero-mean, unit-variance random Gaussian matrix similar to the previous experiment. Then, we use SVD to modify the singular values to obtain a matrix with a desired condition number. This is done by linearly scaling  the differences of all singular values with the smallest non-zero singular value. Table~\ref{tab:res3} and Table~\ref{tab:res4} summarize the result of our experiment. In particular, Table~\ref{tab:res3} shows the number of required iterations, while Table~\ref{tab:res4} presents the required time. We observe that as the condition number grows, the number of required iterations and CPU time increases. DSGs requires a fewer number of iterations compared to other algorithms to reach an error of $10^{-10}$ in almost all scenarios. In terms of CPU time, however, DSGS is not the fastest among algorithms considered.

	\begin{table}
		\centering
		\caption{\textsc{The Average Number of Required Iterations  to Reach the Error Level of $10^{-10}$ For An Overdetermined System with A Certain Given Condition Number.} }
		\begin{tabular}{llllllll}
			\hline \hline
			Size ($m\times n$)&  $k(A)$ &  DSGS &  RK &  RGS &  REK      \\
			\hline
			\hline
			$40\times 20$ & $10^1$& \textbf{9,251} & 16,731 & 16,324 & 19,813  \\
			$40\times 20$ & $10^2$& \textbf{429,510}& 1,521,400  & 1,530,549 & 1,779,000   \\
			$40\times 20$ & $10^3$&\textbf{42,503,000} &147,368,000 & 152,590,000&175,640,000 \\
			$40\times 30$   & $10^1$&\textbf{18,457} &24,742  & 24,668 &  29,461   \\
			$40\times 30$   & $10^1$& \textbf{685,800} &2,216,400  & 2,217,500 & 2,615,400    \\
			$40\times 30$   & $10^3$&\textbf{65,444,000}&  220,290,000&217,826,000  &  263,440,000   \\
			$60\times 30$  & $10^1$& \textbf{18,512}& 24,338 & 24,569 & 28,206  \\
			$60\times 30$  & $10^2$&\textbf{756,110}& 2,206,800 & 2,262,800 & 2,708,700  \\
			$60\times 30$  & $10^3$&\textbf{66,628,000}& 212,990,000 &220,510,000 &  256,190,000 \\
			$60\times 50$ & $10^1$&  \textbf{39,538}& 40,330& 40,930 & 47,523 \\
			$60\times 50$ & $10^2$&\textbf{1,362,900}&3,644,200 & 3,797,800 &  4,407,700\\
			$60\times 50$ & $10^3$&\textbf{89,837,000} & 363,740,000& 366,690,000 & 417,850,000 \\
			$80\times 40$ & $10^1$&\textbf{28,435} & 31,764 & 32,177  &  38,858 \\
			$80\times 40$ & $10^2$& \textbf{911,460}&2,989,500 & 2,927,400 &  3,486,000 \\
			$80\times 40$ & $10^3$&\textbf{72,912,000}& 296,892,000 & 298,320,000  & 359,853,000  \\
			$80\times 60$  & $10^1$& 54,069 &49,828 &\textbf{47,672} & 55,977 \\
			$80\times 60$  & $10^2$&\textbf{1,509,100} &4,284,300 &4,469,300 &5,270,600\\
			$80\times 60$  & $10^3$&\textbf{120,660,000} &437,880,000 & 452,480,000 &530,308,000  \\
			\hline
		\end{tabular}
		\label{tab:res3}
		\vspace{1ex}
	\end{table}

	\begin{table}
	\centering
	\caption{\textsc{The Average Required CPU Time in Seconds to Reach the Error Level of $10^{-10}$ For An Overdetermined System.} }
	\begin{tabular}{llllllll}
		\hline \hline
		Size ($m\times n$)& $k(A)$ &   DSGS &  RK &  RGS &  REK      \\
		\hline
		\hline
		$40\times 20$ & $10^1$
		& 0.3750 &0.9175  & 1.0575 &  \textbf{0.0875} \\
		$40\times 20$ & $10^2$
		&17.9780 &  \textbf{7.44}&  13.83&   7.9240\\
		$40\times 20$ & $10^3$
		& 1,517.7& \textbf{570.832}& 1,154.8&798.3517 \\
		$40\times 30$   & $10^1$
		&1.5960 &  1.0190&1.2033  &   \textbf{0.1880}  \\
		$40\times 30$   & $10^2$
		&59.6680 & \textbf{12.9630} & 26.8070 &  16.1760   \\
		$40\times 30$   & $10^3$
		& 5,618 & \textbf{1,216.6} &2,393.3  &  1,648.9   \\
		$60\times 30$  & $10^1$
		&2.3450 & 1.0780 & 1.2530 & \textbf{0.2520}  \\
		$60\times 30$  & $10^2$
		&93.8430 & \textbf{19.0680} & 35.4230 & 23.9490  \\
		$60\times 30$  & $10^3$
		& 8,493.8& 1,986.6 &3,691.2  & \textbf{1,986.6}  \\
		$60\times 50$ & $10^1$
		& 12.6430&1.4200 & 2.0288 & \textbf{0.6687} \\
		$60\times 50$ & $10^2$
		& 463.2130&\textbf{49.4850} & 107.6980 & 61.9989 \\
		$60\times 50$ & $10^3$
		& 26,367& \textbf{5,026.4}& 10,146 & 6,191 \\
		$80\times 40$ & $10^1$
		&9.0680 & 1.3422 &  1.6833 &  \textbf{0.5980} \\
		$80\times 40$ & $10^2$
		&270.1660 & \textbf{44.0870} & 93.7690  & 53.4530  \\
		$80\times 40$ & $10^3$
		& 21,890&  \textbf{4,528.5}&  9,194.9 &  5,559 \\
		$80\times 60$  & $10^1$
		& 3.4236& 1.8625&2.9329  & \textbf{1.2114} \\
		$80\times 60$  & $10^2$
		& 966.9200& \textbf{88.9450}&  196.6740& 117.0433\\
		$80\times 60$  & $10^3$
		& 77,498 & \textbf{9,090.2} & 19,525 & 12,155  \\
		\hline
	\end{tabular}
	\label{tab:res4}
	\vspace{1ex}
\end{table}

In the second part of our experiments, we consider under-determined systems with $m\le n$. Similar to the previous case, we first generate the entries of $A$ independently according to a zero-mean, unit-variance random Gaussian distribution. Here, we do not include S2CD since it only works for strongly convex problems. The termination rule in this case is chosen as $\|Ax-b\|\leq 10^{-10}$. 
It is observed from Table~\ref{tab:res2} that for all  tested settings, the proposed DSGS algorithm outperforms other algorithms. Similar to the previous case, we compare the performance of different algorithms for different  conditions numbers in Fig.~\ref{tab:res5} and Fig.~\ref{tab:res6}.

		\begin{table}
			\centering
			\caption{\textsc{The Average  Number of Required Iterations  to Reach An Error of $10^{-10}$ For An Underdetermined System.}}
			\begin{tabular}{llllll}
				\hline \hline
				Size  ($m\times n$) & DSGS       & RK & RGS & REK       \\
				\hline
				\hline
				$40\times 60$  & {\bf 11,354} & 23,257 & 23,097 & 26,506\\
				$40\times 100$    & {\bf 3,141} & 6,297 & 6,267 & 6,971   \\
				$80\times 100$  & {\bf 65,030} & 124,724 & 128,957 & 144,894\\
				$80\times 140$  & {\bf 15,171} & 30,664 & 30,168 & 34,608\\
				$160\times 200$ & {\bf 151,858} & 274,537 & 281,331 & 328,680 \\
				$160\times 300$  & {\bf 25,692} & 49,791 & 49,908 & 56,438\\
				\hline
			\end{tabular}
			\label{tab:res2}
			\vspace{1ex}
		\end{table}

\begin{table}
	\centering
	\caption{\textsc{The Average Number of Required Iterations  to Reach An Error of $10^{-10}$ For An Underdetermined System with A Certain Given Condition Number.} }
	\begin{tabular}{llllllll}
		\hline \hline
		Size ($m\times n$)&  $k(A)$ &  DSGS &  RK &  RGS &  REK      \\
		\hline
		\hline
		$40\times 60$ & $10^1$
		& \textbf{15,011} & 32,261 & 32,901& 37,957   \\
		$40\times 60$ & $10^2$
		& \textbf{1,373,800} & 3,304,000 & 3,468,600 & 3,999,500  \\
		$40\times 60$ & $10^3$
		& \textbf{137,686,180}& 305,210,000& 306,252,290,&354,600,000 \\
		$80\times 140$   & $10^1$
		&\textbf{29,600} & 65,831 & 66,201 &  76,367   \\
		$80\times 140$   & $10^2$
		&\textbf{2,829,840} &5,465,486  &6,018,850  &7,000,930     \\
		$80\times 140$   & $10^3$
		& \textbf{249,641,600} & 549,779,370 &  572,884,321& 726,486,714    \\
		\hline
	\end{tabular}
	\label{tab:res5}
	\vspace{1ex}
\end{table}

\begin{table}
	\centering
	\caption{\textsc{The Average Required CPU Time in Seconds to Reach An Error of $10^{-10}$ For An Underdetermined System.} }
	\begin{tabular}{llllllll}
		\hline \hline
		Size ($m\times n$)& $k(A)$ &   DSGS &  RK &  RGS &  REK      \\
		\hline
		\hline
		$40\times 60$ & $10^1$
		& 3.72 & 1.1360 & 1.6080 &\textbf{0.4620}    \\
		$40\times 60$ & $10^2$
		&\textbf{42.9820} & 52.9140 & 120.4680 &  61.6675 \\
		$40\times 60$ & $10^3$
		& 41,824& \textbf{3,195}& 6,589& 3,971\\
		$80\times 140$   & $10^1$
		&\textbf{3.4400} & 4.1980 & 8.8 & 4.1400    \\
		$80\times 140$   & $10^2$
		&332.2200 & \textbf{158.7440} & 368.0673 &  208.7914   \\
		$80\times 140$   & $10^3$
		&  27,558& \textbf{23,871} & 47,411 &  31,097   \\
		\hline
	\end{tabular}
	\label{tab:res6}
	\vspace{1ex}
\end{table}

\begin{figure}[htb]
			\centering
			\includegraphics[width=8cm]{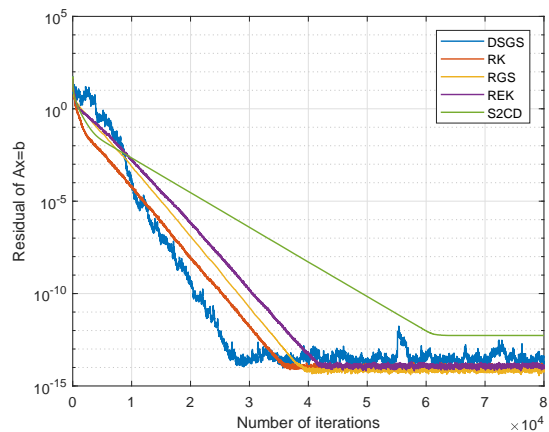}
			\caption{Convergence behavior of different algorithms for an underdetermined system of size $40\times 100$.} \label{fig:sam}
		\end{figure}

\begin{figure}[htb]
			\centering
			\includegraphics[width=1\textwidth]{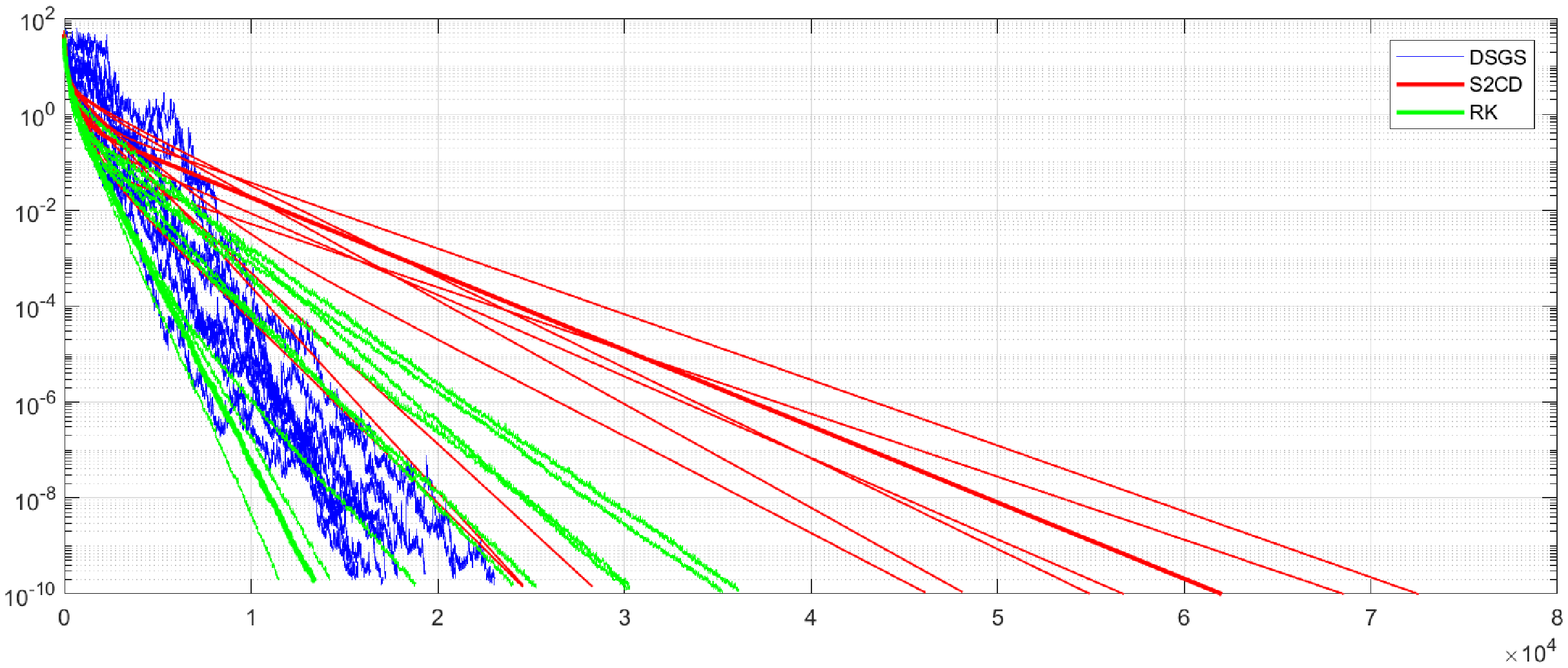}
			\caption{Convergence behavior for DSGS, RK and RGS algorithms for a linear system of size $40\times 30$. In this figure, $10$ realizations are shown for each of the algorithms.}\label{fig:over1}
			\bigbreak
			\includegraphics[width=1\textwidth]{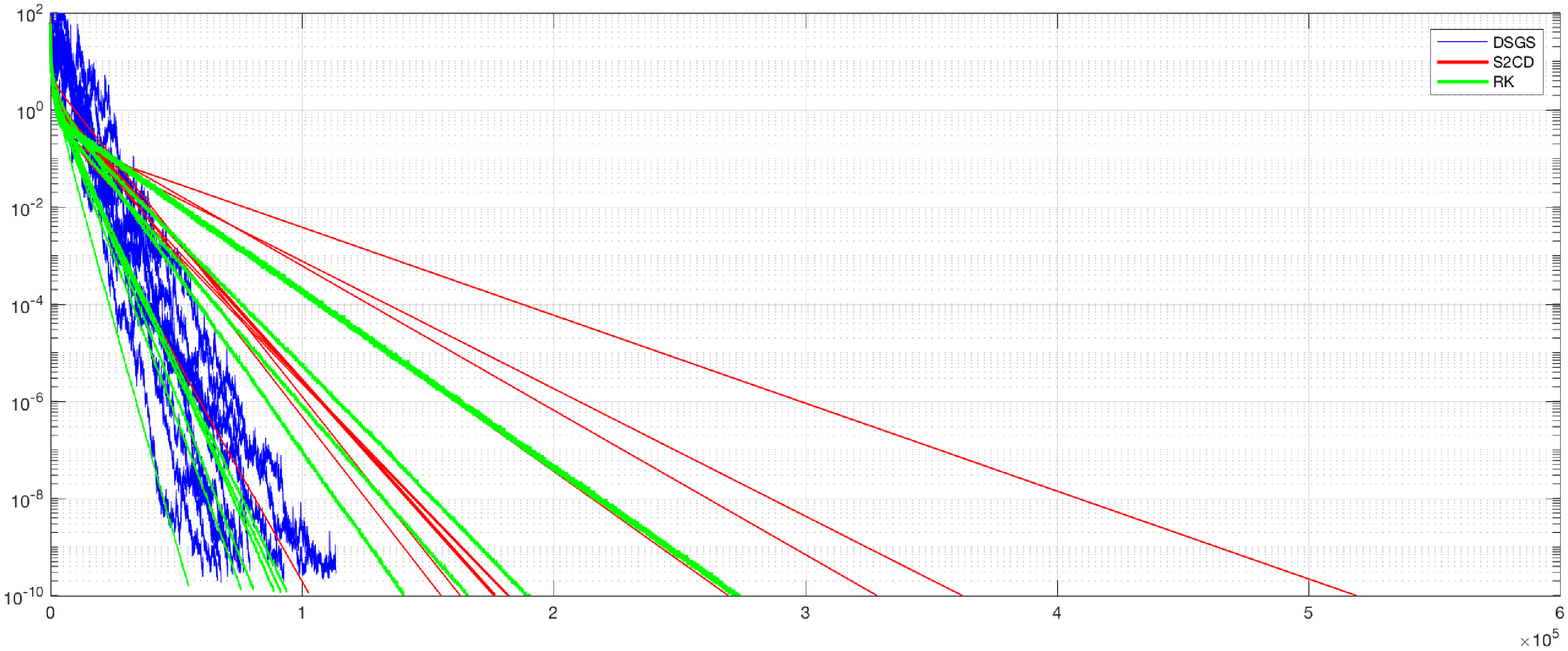}
			\caption{Convergence rate for DSGS, RK and RGS algorithms for a linear system of size $60\times 50$. In this figure, $10$ realizations are shown for each of the algorithms.}\label{fig:over2}
		\end{figure}

		We have also illustrated the convergence of various algorithms in terms of least squares error versus iteration number in Fig.~\ref{fig:sam}. This plot is generated for a random Gaussian matrix of size $40\times30$. It is observed that for the considered linear system with random matrix, DSGS slightly outperforms other algorithms. Since DSGS chooses each row and coordinate to update randomly, the residual of least squares in DSGS algorithm does not decrease monotonically.   We have repeated this experiment in  Fig.~\ref{fig:over1} and Fig.~\ref{fig:over2} over ten different realization of matrix~$A$. The linear system in Fig. \ref{fig:over1} is of size $40 \times 30$ and it is seen that DSGS and RK have close convergence behavior and both faster than that  of S2CD. We increase the size of  matrix $A$ to $60\times 50$ and similar observation has been made in Fig.~\ref{fig:over2}.
		
\begin{figure}[htb]
			\centering
			\begin{minipage}{.5\textwidth}
				\includegraphics[width=1\linewidth]{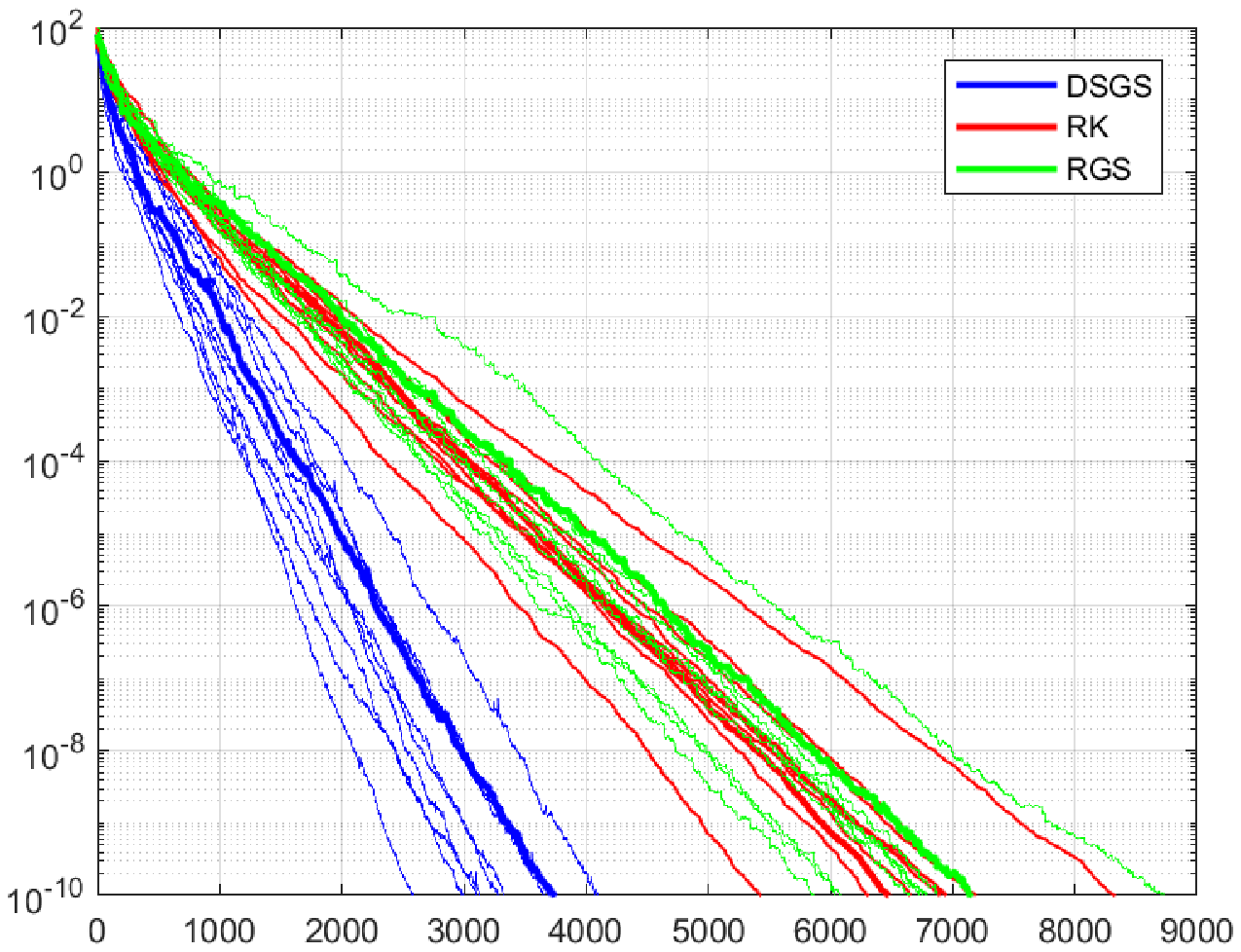}
				\end{minipage}~
			\begin{minipage}{.5\textwidth}
				\includegraphics[width=1\linewidth]{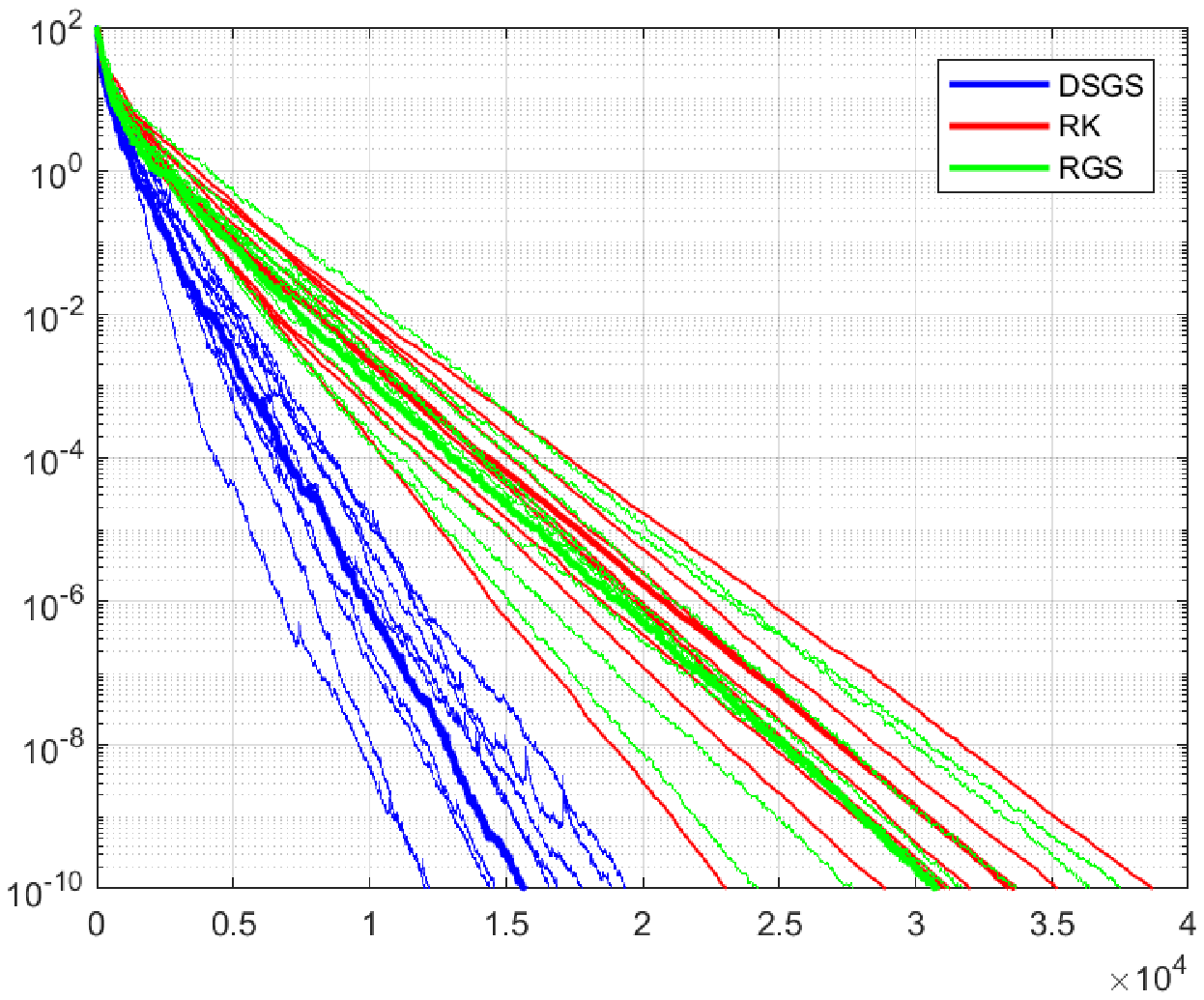}
				\end{minipage}
							\caption{Convergence rate for DSGS, RK and RGS algorithms for a linear system with a matrix of size $40\times 100$ and $80\times 140$.}
							\label{fig:under}
		\end{figure}

		Finally, in Fig. \ref{fig:sam} we compare the convergence behavior of different algorithms for an under-determined system of size $40\times 100$. Similar to the results given in Table \ref{tab:res2}, it is seen that DSGD requires smaller number of iterations to converge. We also plot a number of  different cases with $10$ realizations for each of the algorithms in Figs. \ref{fig:over1} -- \ref{fig:under}.  We observe that  S2CD is only designed for strongly convex problems, therefore it does not perform well for the under-determined systems whose the objective function is  convex but not strongly convex. 

		
	\bibliographystyle{unsrt}

\end{document}